\documentclass[12pt]{article}
\usepackage{amsthm,amsfonts,amssymb,amscd}

\begin{document}

\begin{center}
\Large \bf Birationally rigid varieties \\
with a pencil of Fano double covers. III
\end{center}
\vspace{1cm}

\centerline{A.V. Pukhlikov %\footnote{Работа выполнена при поддержке
%гранта МД-328.2003.01 для молодых российских ученых-докторов наук
%и гранта Фонда содействия отечественной науке}
}

\parshape=1
3cm 10cm \noindent {\small \quad \quad \quad
\quad\quad\quad\quad\quad\quad\quad {\bf }\newline We complete
the study of birational geometry of Fano fiber spaces $\pi\colon
V\to {\mathbb P}^1$, the fiber of which is a Fano double
hypersur\-face of index 1. For each family of these varieties we
either prove birational rigidity or produce explicitly
non-trivial struc\-tures of Fano fiber spaces. A new linear method
of studying movable systems on Fano fiber spaces $V/{\mathbb P}^1$
is developed. \par Bibliography: 17 titles.} \vspace{1cm}

\section*{Introduction}

The present paper is the concluding one in the series of papers
investigating birational geometry of fiber spaces $\pi\colon
V\to{\mathbb P}^1$, the fiber of which is a double Fano
hypersurface of index 1. For almost all families (except for a
finite set) birational rigidity was proved in [1,2]. The remaining
families are considered in this paper.\vspace{0.3cm}

{\bf 0.1. The list of varieties under consideration.} Recall [2]
that the parameters of a family of fiber spaces $V/{\mathbb P}^1$
are written down in the following format:
$$
((a_1,\dots,a_{M+1}),(a_Q,a_W))
$$
(respectively, $((a_1,\dots,a_M),a_W)$, if the pencils of Fano
double spaces are considered), where
$$
\sigma\colon V\to Q\subset X={\mathbb P}({\cal E})
$$
is a double cover, ${\cal E}=\oplus{\cal O}_{{\mathbb P}^1}(a_i)$,
$a_0=0\leq a_1\leq\dots\leq a_{M+1}$, the morphism $\sigma$ is
ramified over the smooth divisor $W_Q=W\cap Q$,
$$
Q\sim mL_X+a_QR,\,\,W\sim 2lL_X+2a_WR,
$$
$a_Q,a_W\in{\mathbb Z}_+$, $R$ is the class of a fiber of the
projection $\pi_X\colon X\to{\mathbb P}^1$, whereas $L_X\in
\mathop{\rm Pic}X$ is the class of the tautological sheaf of the
fibration $X/{\mathbb P}^1$. For compactness of notations in the
set $(a_1,\dots,a_{M+1})$ we write down only non-zero entries, if
there are any, otherwise we write $(0)$. Here is the list of
families that were excluded from consideration in [1,2]:
$$
\begin{array}{cl}
1. & ((0),(1,0)); \\
2. & ((0),(0,1)); \\
3. & ((0),(0,0)); \\
4. & ((1),(1,0)); \\
5. & ((1),(0,0)); \\
6. & ((1,1),(1,0)); \\
7. & ((1,1),(0,0)).
\end{array}
$$
Besides, in [1,2] the following families of varieties with a
pencil of double spaces were excluded from consideration:
$$
\begin{array}{cl}
1^*. &  ((0),1); \\
2^*. & ((0),0); \\
3^*. & ((1),0); \\
4^*. & ((1,1),0).
\end{array}
$$
In Sec. 1 of this paper we prove birational superrigidity of
general varieties in the families $6,7$ and $4^*$. For varieties
in the families $1-3,5$ and $1^*-3^*$ we construct infinitely many
non-trivial structures of Fano fiber spaces, which excludes
birational rigidity. This solves the problem of birational
rigidity for all families except for the family 4. Those
varieties will be considered in a separate paper. \vspace{0.3cm}

{\bf 0.2. Further results.} To prove certain crucial facts in
[1,2] we used the modern techniques to the very limit. This is
true, in the first place, for the proof of Proposition 3.1 in
[1]. Thus it is natural to try to improve certain technical parts
of the arguments. In Sec. 2 we discuss how the local inequality,
on which the proof of Proposition 3.1 in [1] is based, could be
improved. A stronger version of that inequality is formulated as
a conjecture and completely proved in the case when the graph of
the corresponding sequence of blow ups is a chain.

In Sec. 3 we develop a new method of proving birational rigidity
of Fano fiber spaces $V/{\mathbb P}^1$ (more exactly, a new
method of studying movable linear systems on these fiber spaces),
that is, the {\it linear method}. The name emphasizes its being
different from the traditional {\it quadratic} method, based on
the operation of taking the self-intersection of a movable linear
system $\Sigma$, that is, the operation of making the effective
algebraic cycle $Z=(D_1\circ D_2)$ of codimension 2, where
$D_1,D_2\in\Sigma$ are general divisors. Naturally, all
constructions involved in the quadratic method are quadratic in
the parameters of the system $\Sigma$ (e.g. the $4n^2$-inequality
and other crucial facts). On the contrary, the linear method is
based on the operation of restricting the movable linear system
$\Sigma$ on a certain fiber $F=F_t=\pi^{-1}(t)$ of the fiber
space $V/{\mathbb P}^1$. It is clear that the principal
computations of the linear method are linear in the parameters of
the system $\Sigma$ (degrees, multiplicities, etc.). The linear
method was for the first time successfully applied in [3]. Its
idea could be illustrated in the following way.

{\bf Definition 0.1. (see [3]).} We say that a primitive Fano
variety $F$ is {\it divisorially canonical}, or satisfies the
condition ($C$) (respectively, is {\it divisorially log
canonical}, or satisfies the condition ($L$)), if for any
effective divisor $D\in|-nK_F|$, $n\geq 1$, the pair
\begin{equation}\label{i2}
(F,\frac{1}{n}D)
\end{equation}
has canonical (respectively, log canonical) singularities. If the
pair (\ref{i2}) has canonical singularities for a general divisor
$D\in \Sigma \subset|-nK_F|$ of any {\it movable} linear system
$\Sigma$, then we say that $F$ satisfies the condition of {\it
movable canonicity}, or the condition ($M$).

Let $V/{\mathbb P}^1$ be a Fano fiber space, $\Sigma\subset
|-nK_V+lF|$ a movable linear system on $V$, which is not made
from the pencil of fibers of the projection $\pi$, that is,
$n\geq 1$. Assume that the log pair  $(V,\frac{1}{n}\Sigma)$ is
not canonical, that is, for some exceptional divisor
$E\subset\widetilde V$, where $\varphi\colon \widetilde V\to V$
is a sequence of blow ups, the Noether-Fano inequality
$$
 \nu_E(\Sigma)>na(E)
$$
is satisfied. Assume that the centre $B=\varphi(E)\subset V$ lies
in some fiber $F=F_t$ and is of codimension $\geq 3$ on $V$ (this
case is the hardest one in all particular problems). Since the
linear system $\Sigma$ is movable, its restriction
$\Sigma_F=\Sigma\,|\,_F$ is a non-empty linear system (which may
have fixed components). According to the inversion of adjunction
[4], the log pair
$$
(F,\frac{1}{n}\Sigma_F)
$$
is not log canonical. Suppose that we knew from the start that
the fibers $F$ satisfy the condition ($L$) (the more so, the
condition ($C$)). Now we get a contradiction excluding the maximal
singularity $E$.

This is the general outline of the scheme of the linear method.
In order to realize it, one needs to be able to prove the
condition ($L$), which in most cases is very hard. The only known
way to do it was shown in [3]. It requires essentially stronger
conditions of general position (regularity) than those used by
the quadratic method [5]. However, the linear method has certain
advantages compared to the quadratic one, namely, it absolutely
does not use the condition of ``twistedness'' of the fiber space
$V/{\mathbb P}^1$ over the base. The only information which is
used is the properties of fibers.

In Sec. 3 by means of the linear method we prove birational
superrigidity of the fiber spaces $V/{\mathbb P}^1$ into double
spaces, double quadrics and Fano hypersurfaces of index 1.
\vspace{0.3cm}

{\bf 0.3. Acknowledgements.} A considerable part of results of
the present paper (in particular, the computations in Sec. 2 and
the key parts of the linear method) was obtained by the author
during his stay at Max-Planck-Institut f\" ur Mathematik in Bonn
in the summer and autumn of 2003. The author is grateful to the
institute for the hospitality and the excellent conditions of
work.

%%%%%%%%%%%%%%%%%%%%%%%%%%%%%%%%%%%%%%%%%%%%%%%%%%%%%%%%%%%%%%%
%%%%%%%%%%%%%%%%%%%%%%%%%%%%%%%%%%%%%%%%%%%%%%%%%%%%%%%%%%%%%%%
%%%%%%%%%%%%%%%%% section 1
\section{Birationally rigid and non-rigid Fano fiber spaces}

{\bf 1.1. Varieties of the type ${\bf ((1,1),(1,0))}$.} Here
$X={\mathbb P}({\cal E})$, where the sheaf ${\cal E}$ is of the
form ${\cal E}={\cal O}^{\oplus M}_{{\mathbb P}^1}\oplus{\cal
O}_{{\mathbb P}^1}(1)^{\oplus 2}$. The space $H^0(X,{\cal
L}_X\otimes\pi^*{\cal O}_{{\mathbb P}^1}(-1))$ is two-dimensional
and defines a pencil of divisors $|L_X-R|$. Its base set
$\Delta_X=\mathop{\rm Bs}|L_X-R|$ is of codimension 2: it is easy
to see that
$$
\Delta_X={\mathbb P}({\cal O}^{\oplus M}_{{\mathbb
P}^1})\cong{\mathbb P}^{M-1}\times{\mathbb P}^1.
$$
Furthermore, $Q\sim mL_X+R$ and $W\sim 2lL_X$. Set
$$
\Delta_Q=\Delta_X\cap Q,\quad \Delta=\sigma^{-1}(\Delta_Q)\subset
V.
$$
Obviously, $\Delta_Q$ is a smooth divisor of bidegree $(m,1)$ on
$\Delta_x={\mathbb P}^{M-1}\times{\mathbb P}^1$, $\Delta\subset V$
is a smooth irreducible subvariety of codimension 2.

{\bf Lemma 1.1.} {\it The anticanonical linear system $|-K_V|$ is
movable, and moreover $\mathop{\rm Bs}|-K_V|=\Delta$. Furthermore,
$$
-K_V\in\partial A^1_{\mathop{\rm mov}}V.
$$
More precisely, every linear system $|-nK_V+lF|$ is empty for
$l<0$.}

{\bf Proof.} The first claim of the lemma follows from the fact
that the anticanonical linear system $|-K_V|=|L_V-F|$.

Assume that the linear system $\Sigma=|-nK_V+lF|$ is non-empty
for some $l<0$. We must show that this assumption leads to a
contradiction. In order to do it, we construct a special family
of surfaces sweeping out $V$. Restricting the linear system
$\Sigma$ onto a general surface of this family, we must get a
non-empty linear system of curves. As we will see, the latter is
impossible.

Let
$$
\alpha\colon{\cal O}^{\oplus M}_{{\mathbb P}^1}\to{\cal
O}_{{\mathbb P}^1}\oplus{\cal O}_{{\mathbb P}^1}\quad
\mbox{and}\quad\beta\colon{\cal O}_{{\mathbb P}^1}(1)^{\oplus
2}\to{\cal O}_{{\mathbb P}^1}(1)
$$
be surjective morphisms of sheaves. Their direct sum
$\alpha\oplus\beta$ determines an inclusion of projective bundles
$$
S=S(\alpha,\beta)={\mathbb P}({\cal O}^{\oplus 2}_{{\mathbb
P}^1}\oplus{\cal O}_{{\mathbb P}^1}(1))\hookrightarrow X={\mathbb
P} ({\cal E}).
$$
The subvariety $S$ is a ${\mathbb P}^2$-bundle over ${\mathbb
P}^1$, where the intersection $S\cap\Delta_X\cong {\mathbb
P}^1\times{\mathbb P}^1$ is a {\it divisor} on $S$.

Set
$$
V_S=\sigma^{-1}(Q\cap S),\quad\Delta_S=\Delta\cap V_S.
$$
For sufficiently general $\alpha,\beta$ the surface $V_S$ is
smooth, $\Delta_S\subset V_S$ is a smooth curve. Moreover, the
surfaces $V_S$ sweep out the variety $V$.

Obviously, $\sigma(\Delta_S)=\Delta_Q\cap S\subset{\mathbb
P}^1\times{\mathbb P}^1$ is a curve of bidegree $(m,1)$. Denote
by the symbols $L_S$ and $F_S$ the restrictions of the classes
$L_V$ and $F$ onto the surface $V_S$, respectively.

{\bf Lemma 1.2.} (i) {\it The following equivalences hold:}
$$
\Delta_S\sim L_S-F_S\sim(-K_V)|_{V_S}.
$$

(ii) {\it The following equality holds:
$(\Delta_S\cdot\Delta_S)=2(1-m)$.}

{\bf Proof.} (i) By the construction of the variety $S$ we get:
$$
S\cap\Delta_X\sim L_X-R,
$$
the rest is obvious.

(ii) By the part (i) we get
$$
(\Delta_S\cdot\Delta_S)=(\Delta_S\cdot L_S)-(\Delta_S\cdot F_S).
$$
As we noted above, $\sigma(\Delta_S)$ is a curve of bidegree
$(m,1)$ on $S\cap\Delta_X={\mathbb P}^1\times{\mathbb P}^1$. It
is easy to see that the restriction of the tautological sheaf
$L_X|_{S\cap\Delta_X}$ is of bidegree $(1,0)$. Therefore
$(\Delta_S\cdot L_S)=2$ and $(\Delta_S\cdot F_S)=2m$. (The
intersection indices are doubled because of the double cover
$\sigma$.) Q.E.D. for the lemma.

Let us complete the proof of Lemma 1.1. Let $\Sigma_S$ be the
restriction of the linear system $\Sigma$ onto the surface $V_S$.
This is a non-empty linear system of curves. By the previous
lemma,
$$
\Sigma_S\subset |n\Delta_S+lF_S|,
$$
where $n\geq 0$, $l<0$. Take a general curve $C\in\Sigma_S$ and
write down
$$
C=m^+\Delta_S+C^+,
$$
where $n^+\in{\mathbb Z}_+$ and all irreducible components of the
curve $C^+$ are different from $\Delta_S$. Of course, $m^+\leq n$
(otherwise the class $(-F_S)$ would have been effective).
Consequently,
$$
C^+\sim (n-m^+)\Delta_S+lF_S,
$$
so that
$$
0\leq (C^+\cdot\Delta_S)=2(n-n^+)(1-m)+2ml<0
$$
(recall that $l<0$). This contradiction completes the proof of
Lemma 1.1.

{\bf Lemma 1.3.} {\it The Fano fiber space $V/{\mathbb P}^1$
satisfies the generalized $K^2$-condition of depth $1/m$, that is,
$$
K^2_V-\frac{1}{m}H_F\notin\mathop{\rm Int}A^2_+V,
$$
where $H_F\in A^2_{\mathbb R}V$ is the class of a hyperplane
section of a fiber, that is, $H_F=(-K_V\cdot F)$.}

{\bf Proof.} By the formula (1) of the paper [2], $(K^2_V\cdot
L^{M-1}_V)=2$. Since for the degree of the fiber we have obviously
$(F\cdot L^M_V)=(H_F\cdot L^{M-1}_V)=2m$, we get
$$
((K^2_V-\frac{1}{m}H_F)\cdot L^{M-1}_V)=0
$$
which immediately implies the claim of the lemma.

{\bf Corollary 1.1.} {\it A general Fano fiber space $V/{\mathbb
P}^1$ of type $((1,1),(1,0))$ is birationally superrigid. The
projection $\pi\colon V\to{\mathbb P}^1$ gives the only
non-trivial structure of a rationally connected fiber space on
$V$. The groups of birational and biregular self-maps of the
variety $V$ coincide:}
$$
\mathop{\rm Bir}V=\mathop{\rm Aut}V={\mathbb Z}/2{\mathbb Z}
$$

{\bf Proof:} this follows immediately from Theorem 2 of the paper
[2] ($1/m<2$, so that the fiber space $V/{\mathbb P}^1$ satisfies
the $K^2$-condition of depth 2, whereas any movable linear system
$\Sigma$ is a subsystem of the complete linear system
$|-nK_V+lF|$ with $n,l\in{\mathbb Z}_+$).\vspace{0.3cm}

%%%%%%%%%%%%%%%%%%%%%%%%%%%%%%%%%%%%%%%%%%%%%%%%%%%%%%%%%%%%%%%%%%
%%%%%%%%%%%%%%%%%%%%% subsection 1.2

{\bf 1.2. Varieties of type ((1,1),(0,0)).} Birational geometry
of varieties of this type is somewhat more complicated than
birational geometry of the varieties of type ${((1,1),(1,0))}$.
The projective bundle $X$ and the locally free sheaf ${\cal E}$
are the same as in Sec. 1.1. However, in the case under
consideration $Q\sim mL_X$, so that $-K_V=L_V$ and thus the
linear system
$$
|-K_V-F|=\sigma^*(|L_X-R|\,\Bigr|_Q)
$$
is movable. Consequently, the variety $V$ does not satisfy the
$K$-condition. Let
$$
\varphi\colon V\dashrightarrow{\mathbb P}^1
$$
be the rational map, given by the pencil $|-K_V-F|$. Birational
geometry of the variety $V$ is completely described by

{\bf Proposition 1.1.} (i) {\it The variety $V$ is birationally
superrigid: for any movable linear system $\Sigma$ on $V$ its
virtual and actual thresholds of canonical adjunction coincide,}
$$
c_{\mathop{\rm virt}}(\Sigma)=c(\Sigma).
$$

(ii) {\it On the variety $V$ there are exactly two non-trivial
structures of a rationally connected fiber space, namely
$\pi\colon V\to{\mathbb P}^1$ and $\varphi\colon
V\dashrightarrow{\mathbb P}^1$. These structures are birationally
distinct, that is, there is no birational self-map
$\chi\in\mathop{\rm Bir}V$, transforming the fibers of $\pi$ into
the fibers of $\varphi$. The groups of birational and biregular
self-maps of the variety $V$ coincide:} $\mathop{\rm
Bir}V=\mathop{\rm Aut}V$.

(iii) {\it There is a unique, up to a fiber-wise isomorphism,
Fano fiber space $\pi^+\colon V^+\to{\mathbb P}^1$ of the same
type $((1,1),(0,0))$, such that the following diagram commutes:
$$
\begin{array}{ccccc}
& V &  \stackrel{\chi}{\dashrightarrow} & V^+ & \\
\varphi &  \downarrow  &    & \downarrow &  \pi^+ \\
   & {\mathbb P}^1 & = & {\mathbb P}^1, &
\end{array}
$$
where $\chi$ is a birational map. The correspondence $V\to V^+$
is involutive, that is,} $(V^+)^+=V$.

{\bf Proof.} The symbols $\Delta_Q$, $\Delta_X$ and $\Delta$ mean
the same as above (Sec. 1.1). However, in the case under
consideration $\Delta_Q$ is a smooth divisor of bidegree $(m,0)$
on $\Delta_X={\mathbb P}^{M-1}\times{\mathbb P}^1$, that is,
$\Delta_Q=\Delta_{\mathbb P}\times{\mathbb P}^1$, where
$\Delta_{\mathbb P}\subset{\mathbb P}^{M-1}$ is a smooth
hypersurface of degree $m$. Similarly,
$\Delta=\Delta_F\times{\mathbb P}^1$, where
$\Delta_F=\sigma^{-1}(\Delta_{\mathbb P})$ is the double cover of
the hypersurface $\Delta_{\mathbb P}$, branched over
$W\cap\Delta_{\mathbb P}$.

{\bf Lemma 1.4.} {\it The base set of the movable linear system
$|-K_V-F|$ is $\mathop{\rm Bs}|-K_V-F|=\Delta$. Furthermore,
$$
-K_V-F\in\partial A^1_{\mathop{\rm mov}}V.
$$
More precisely, $|-nK_V+lF|=\emptyset$ for $l<-n$.}

{\bf Proof} is completely similar to the proof of Lemma 1.1. The
only difference is that this time the curve $\Delta_S\sim L_S-F_S$
(the surface $S$ is constructed in absolutely the same way as in
the proof of Lemma 1.1) is not irreducible. The curve $\Delta_S$
is a union of $2m$ disjoint $(-1)$-curves,
$$
\Delta_S=\sum^m_{i=1}(C^+_i+C^-_i),\quad\Delta_Q
=\sum^m_{i=1}C_i,\quad\sigma(C^{\pm}_i)=C_i.
$$
The lines $C_i\subset S\cap\Delta_X={\mathbb P}^1\times{\mathbb
P}^1$ are of bidegree $(1,0)$, they correspond to the $m$ points
of intersection $\Delta_{\mathbb P}\cap S$. Note that the branch
divisor $W\cap\Delta_X\subset{\mathbb P}^{M-1}\times{\mathbb
P}^1$ is of bidegree $(2l,0)$, that is,
$$
W\cap\Delta_X=W_{\mathbb P}\times{\mathbb P}^1,
$$
where $W_{\mathbb P}\subset{\mathbb P}^{M-1}$ is a general
hypersurface of degree $2l$. Thus for a general surface $S$ the
lines $C_1,\dots,C_m$ lie outside the branch divisor, so that we
get $C^+_i\cap C^-_i=\emptyset$.

Now we prove that the linear system $|aL_V-bF|$ is empty for
$b>a$ in exactly the same way as in the proof of Lemma 1.1 (where
the irreducible curve $\Delta_S$ is replaced by $2m$ disjoint
$(-1)$-curves $C^{\pm}_i$).

Proof of Lemma 1.4 is complete.

Now let us study the rational map $\varphi\colon V
\dashrightarrow {\mathbb P}^1$. In order to do that, we need an
explicit coordinate presentation of the varieties $X$, $Q$ and
$W$, participating in the construction of the Fano fiber space
$V/{\mathbb P}^1$.

Consider the locally free subsheaves
$$
{\cal E}_0={\cal O}^{\oplus M}_{{\mathbb P}^1}\hookrightarrow{\cal
E}\quad \mbox{and}\quad {\cal E}_1={\cal O}_{{\mathbb
P}^1}(1)^{\oplus 2}\hookrightarrow{\cal E}.
$$
Obviously, ${\cal E}={\cal E}_0\oplus{\cal E}_1$. Let
$\Pi_0\subset H^0(X,{\cal L}_X)$ be the subspace, corresponding
to the space of sections of the sheaf $H^0({\mathbb P}^1,{\cal
E}_0)\hookrightarrow H^0({\mathbb P}^1,{\cal E})$. Set also
$$
\Pi_1=H^0(X,{\cal L}_X\otimes\pi^*{\cal O}_{{\mathbb
P}^1}(-1))=H^0({\mathbb P}^1,{\cal E}_1(-1)).
$$
Let $x_0,\dots,x_{M-1}$ be a basis of the space $\Pi_0$,
$y_0,y_1$ a basis of the space $\Pi_1$. Then the sections
\begin{equation}\label{a1}
x_0,\dots,x_{M-1},y_0t_0,y_0t_1,y_1t_0,y_1t_1,
\end{equation}
where $t_0,t_1$ is a system of homogeneous coordinates on
${\mathbb P}^1$, make a basis of the space $H^0(X,{\cal L}_X)$.
It is easy to see that the complete linear system (\ref{a1})
defines a morphism
$$
\xi\colon X\to\bar X\subset{\mathbb P}^{M+3},
$$
the image $X$ of which is a quadratic cone with the vertex space
${\mathbb P}^{M-1}=\xi(\Delta_X)$ and a smooth quadric in
${\mathbb P}^3$, isomorphic to ${\mathbb P}^1\times{\mathbb
P}^1$, as a base. The morphism $\xi$ is birational, more
precisely,
$$
\xi\colon X\setminus \Delta_X\to\bar X\setminus \xi(\Delta_X)
$$
is an isomorphism and $\xi$ contracts $\Delta_X={\mathbb
P}^{M-1}\times{\mathbb P}^1$ onto the vertex space of the cone.
Let
$$
u_0,\dots,u_{M-1},u_{00},u_{01},u_{10},u_{11}
$$
be the homogeneous coordinates on ${\mathbb P}^{M+3}$,
corresponding to the ordered set of sections (\ref{a1}). The cone
$\bar X$ is given by the equation
$$
u_{00}u_{11}=u_{01}u_{10}.
$$
On the cone $\bar X$ there are two pencils of $(M+1$)-planes,
corresponding to the two pencils of lines on a smooth quadric in
${\mathbb P}^3$. Let $\tau\in\mathop{\rm Aut}{\mathbb P}^{M+3}$
be the automorphism permuting the coordinates $u_{01}$ and
$u_{10}$ and not changing the other coordinates. Obviously,
$\tau\in\mathop{\rm Aut}\bar X$ is an automorphism of the cone
$\bar X$, permuting the above-mentioned pencils of
$(M+1)$-planes. One of these pencils is the image of the pencil
of fibers of the projection $\pi$, that is, the pencil
$\xi(|R|)$. For the other pencil we get the equality
$$
\tau\xi(|R|)=\xi(|L_X-R|).
$$
The automorphism $\tau$ induces an involutive birational self-map
$$
\tau^+\in\mathop{\rm Bir}X.
$$
More precisely, $\tau^+$ is a biregular automorphism outside a
closed subset $\Delta_X$ of codimension 2. Let
$$
\varepsilon\colon\tilde X\to X
$$
be the blow up of the smooth subvariety $\Delta_X$. Obviously,
the variety $\tilde X$ is isomorphic to the blow up of the cone
$\bar X$ at its vertex space $\xi(\Delta_X)$. It is easy to check
that $\tau^+$ extends to a biregular automorphism of the smooth
variety $\tilde X$.

Set $Q^+=\tau^+(Q)\subset X$, $W^+=\tau^+(W)\subset X$. The
divisors $Q^+$ and $W^+$ are well defined because $\tau^+$ is an
isomorphism in codimension 1.

{\bf Lemma 1.5.} {\it The divisors $Q^+$ and $W^+$ are divisors
of general position in the linear systems $|mL_X|$ and $|2lL_X|$,
respectively. In particular, $Q^+$, $W^+$ and $Q^+\cap W^+$ are
smooth varieties.}

{\bf Proof.} The claim follows immediately from the fact that the
linear systems $|kL_X|$, $k\in{\mathbb Z}_+$, are invariant under
$\tau^+$, whereas $Q$ and $W$ are sufficiently general divisors
of the corresponding linear systems. Note that if a divisor
$D\in|kL_X|$ is given by a polynomial
$$
h(u_0,\dots,u_{M-1},u_{00},u_{01},u_{10},u_{11}),
$$
of degree $k$, then its image $\tau^+(D)$ is given by the
polynomial
$$
h^+(u_*)=h(u_0,\dots,u_{M-1},u_{00},u_{10},u_{01},u_{11})
$$
with permuted coordinates $u_{01}$ and $u_{10}$. Q.E.D. for the
lemma.

Let $\sigma^+\colon V^+\to Q^+$ be the double cover, branched
over a smooth divisor $Q^+\cap W^+$. Obviously, $V^+/{\mathbb
P}^1$ is a general Fano fiber space of type $((1,1),(0,0))$.

{\bf Lemma 1.6.} (i) {\it The map $\tau^+$ lifts to a birational
map $\chi\colon V \dashrightarrow V^+$, biregular in codimension
1.}

(ii) {\it The action of $\chi$ on the Picard group is given by
the formulas
$$
\chi^*K_{V^+}=K_V,\quad \chi^*F^+=-K_V-F,
$$
where $F^+$ is the class of the fiber of the projection
$V^+\to{\mathbb P}^1$, so that $\mathop{\rm Pic}V^+={\mathbb
Z}K_{V^+}\oplus{\mathbb Z}F^+$. }

(iii) {\it The construction of the variety $V^+$ is involutive:}
$(V^+)^+\cong V$.

{\bf Proof:} the claims (i)-(iii) are obvious. Just note that the
following presentation holds: $\chi=q^+\circ q^{-1}$, where
$q\colon\tilde V\to V$ and $q^+\colon \tilde V\to V^+$ are blow
ups of the smooth subvarieties of codimension two $\Delta\subset
V$ and $\Delta^+\subset V^+$, respectively. Furthermore,
$E=q^{-1}(\Delta)$ is the exceptional divisor of both blow ups,
$E=\Delta\times{\mathbb P}^1=\Delta_F\times{\mathbb
P}^1\times{\mathbb P}^1$, whereas the projections $q\,|\,_E$ and
$q^+\,|\,_E$ are projections with respect to the second and third
direct factors, respectively.

{\bf Lemma 1.7.} {\it The Fano fiber space $V/{\mathbb P}^1$
satisfies the generalized $K^2$-condition of depth 2.}

{\bf Proof.} This immediately follows from the equality
$$
((K^2_V-2H_F)\cdot L^{M-1}_V)=0.
$$
Q.E.D. for the lemma.

Finally, let us prove Proposition 1.1. Let $\Sigma\subset
|-nK_{V^+}+lF|$ be a movable linear system. If $l\in{\mathbb
Z}_+$, then by Theorem 2 of the paper [2] we get the desired
coincidence of the thresholds: $c_{\mathop{\rm
virt}}(\Sigma)=c(\Sigma)$. Assume that $l<0$. Consider the linear
system $\Sigma^+=\tau^+(\Sigma)$ on $V^+$. By Lemma 1.6,
$\Sigma^+\subset|-n_+K_{V^+}+l_+F^+|$, where
$$
l_+=-l\geq 1.
$$
Since $\tau^+$ is an isomorphism in codimension 1, we get
$c(\Sigma)=c(\Sigma^+)$. Again applying Theorem 2 of the paper
[2], we obtain the desired coincidence of thresholds
$$
c_{\mathop{\rm virt}}(\Sigma^+)=c_{\mathop{\rm
virt}}(\Sigma)=c(\Sigma^+)=c(\Sigma)=n_+=n+l.
$$
This proves birational superrigidity.

Let us prove the claim (ii). The standard argument (see [2, Sec.
1.1]) shows that on $V$ there are exactly two non-trivial
structures of a rationally connected fiber space (the arguments
above imply that if a movable linear system $\Sigma$ satisfies the
equality $c_{\mathop{\rm virt}}(\Sigma)=0$, then either $\Sigma$
is made from the pencil $|F|$, or $\Sigma$ is made from the
pencil $|-K_V-F|$, which gives a description of the existing
structures). For a general variety $V$ these structures cannot be
birationally equivalent. Indeed, by birational superrigidity of
Fano double hypersurfaces of index 1, any birational map
$\chi^+\in\mathop{\rm Bir}V$, which transforms the pencil $|F|$
into the pencil $|-K_V-F|$, induces a biregular isomorphism of
the fibers of general position in the pencils $|F|$ and $|F^+|$
(the latter is taken on the variety $V^+$). Therefore, $\chi^+$
induces a biregular isomorphism of the fibers of general position
of the fiber spaces $Q/{\mathbb P}^1$ and $Q^+/{\mathbb P}^1$.
Now from [6] for $m\geq 3$ we get that these fiber spaces are
globally fiber-wise isomorphic. It checks easily that for a
sufficiently general divisor $Q\subset X$ this is impossible. For
$m=2$ we argue in a similar way, using the branch divisor $W$.

Finally, the claim (iii) follows from the arguments above.

Q.E.D. for Proposition 1.1.\vspace{0.3cm}

%%%%%%%%%%%%%%%%%%%%%%%%%%%%%%%%%%%%%%%%%%%%%%%%%%%%%%%%%%%%%%%%%%%%
%%%%%%%%%%%%%%%%%%%% subsection 1.3.

{\bf 1.3. Varieties of type ((0),(1,0)).} Here $X={\mathbb
P}\times{\mathbb P}^1$, the hypersurface $Q\subset X$ is of
bidegree $(m,1)$, the hypersurface $W\subset X$ is of bidegree
$(2l,0)$, that is, $W=W_{\mathbb P}\times {\mathbb P}^1$, where
$W_{\mathbb P}\subset P$ is a hypersurface of degree $2l$. Let
$(u:v)$ be homogeneous coordinates on ${\mathbb P}^1$, $(x_0:
\dots :x_{M+1})$ be homogeneous coordinates on ${\mathbb P}$. The
hypersurface $Q$ is given by the equation
$$
uf_++vf_-=0,
$$
where $f_{\pm}(x_0,\dots, x_{M+1})$ are homogeneous polynomials
of degree $m$, and the hypersur\-face $W$ is given by the equation
$h(x_*)=0$, $\mathop{\rm deg} h=2l$.

Let $\beta_{\mathbb P}\colon Y_{\mathbb P}\to {\mathbb P}$
(respectively, $\beta_X\colon Y_X\to X$) be the double space,
branched over $W_{\mathbb P}$ (respectively, the double cover,
branched over $W$). Obviously, the direct product $X={\mathbb
P}\times{\mathbb P}^1$ generates the direct product
$Y_X=Y_{\mathbb P}\times {\mathbb P}^1$, compatible with the
double covers $\beta_{\mathbb P}$, $\beta_X$.

{\bf Proposition 1.2.} {\it The projection $q\colon Y_X\to
Y_{\mathbb P}$ onto the first factor determines the birational
morphism
$$
q_V=q|_V\colon V\to Y_{\mathbb P},
$$
contracting the exceptional divisor $E\subset V$:
$$
E=q^{-1}_V(\{\beta^*_{\mathbb P}f_+=\beta^*_{\mathbb P}f_-=0\}).
$$
The birational morphism $q_V$ transforms the pencil $|F|$ of
fibers of the fiber space $V/{\mathbb P}^1$ into the pencil of
divisors
$$
\{\lambda_+\beta^*_{\mathbb P}f_++\lambda_-\beta^*_{\mathbb
P}f_-=0\,|\,\lambda_{\pm}\in{\mathbb C}\}
$$
on $Y_{\mathbb P}$. The inverse birational map $q^{-1}_V$ is the
blow up of the base set of this pencil.}

{\bf Proof:} this is obvious.

{\bf Remark 1.2.} The variety $Y_{\mathbb P}$ is a Fano variety
of index $m+1$. Any structure of a fiber space ${\mathbb
P}\dashrightarrow S$ into Fano complete intersections  of type
$(b_1,\dots,b_e)$ with $b_1+\dots+b_e\leq m$ generates a
structure of a fiber space $Y_{\mathbb P}\dashrightarrow S$ into
Fano varieties of index $m-b_1-\dots-b_e+1\geq 1$. An example of
such structure is given by the pencil $(q_V)_*|F|$ described in
Proposition 1.2. For this method of constructing non-trivial
structures of a fiber space on $Y_{\mathbb P}$ the highest
dimension of the base corresponds to the case $b_1=\dots=b_m=1$,
that is, if we fiber ${\mathbb P}$ into linear subspaces of
codimension $m$. In this case $\mathop{\rm dim}S=m$. The variety
$Y_{\mathbb P}$ is certainly not birationally rigid.

{\bf Conjecture 1.1.} {\it If the hypersurface $W_{\mathbb
P}\subset {\mathbb P}$ is sufficiently general, then for any
structure of a rationally connected fiber space $Y_{\mathbb
P}\dashrightarrow S$ the inequality $\mathop{\rm dim} S\leq m$
holds. If, moreover, $\mathop{\rm dim}S=m$, then there is a
linear projection ${\mathbb P}\dashrightarrow {\mathbb P}^m$ and
a birational map $S\dashrightarrow {\mathbb P}^m$ making the
following diagram commutative:}
$$
\begin{array}{ccc}
Y_{\mathbb P}&\stackrel{2:1}{\longrightarrow} & {\mathbb P}\\
\downarrow & &\downarrow\\
S & \dashrightarrow & {\mathbb P}^m.
\end{array}
$$
\vspace{0.3cm}

%%%%%%%%%%%%%%%%%%%%%%%%%%%%%%%%%%%%%%%%%%%%%%%%%%%%%%%%%%%%%
%%%%%%%%%%%%%%%%%%%%% subsection 1.4

{\bf 1.4. Varieties of type ((0),(0,1)).} Here $X={\mathbb
P}\times{\mathbb P}^1$ and for the variety $Q\subset X$ we also
get the direct decomposition $Q=G\times {\mathbb P}^1\subset
{\mathbb P}\times{\mathbb P}^1$, where $G\subset {\mathbb P}$ is
a hypersurface of degree $m$. Let $q\colon Q\to G$ be the
projection onto the second factor, $L_x=q^{-1}(x)\cong{\mathbb
P}^1$ the fiber over an arbitrary point $x\in G$. Set
$C_x=\sigma^{-1}(L_x)\subset V$. Thus $C_x$ is the fiber of the
projection
$$
q_V=q\circ\sigma\colon V\to G.
$$
The double cover $\sigma\colon C_x\to L_x$ is branched over two
points (because $a_W=1$): the equation of the hypersurface
$W\subset X$ is of the form
$$
f_{uu}u^2+2f_{uv}uv+f_{vv}v^2=0,
$$
where $(u:v)$ are homogeneous coordinates on ${\mathbb P}^1$,
$f_*$ are homogeneous polynomials of degree $2l$ on ${\mathbb
P}$. If the curve $C_x$ is irreducible, then it is a smooth
conic; otherwise, $C_x$ is a pair of smooth rational curves,
provided that all three polynomials $f_{uu}$, $f_{uv}$, $f_{vv}$
do not vanish identically at the point $x$. What has been just
said implies

{\bf Proposition 1.3.} {\it The projection $q_V$ realizes $V$ as
a conic bundle over the base $G\subset{\mathbb P}$. The
discriminant divisor is given by the equation
$f^2_{uv}=f_{uu}f_{vv}$.}

{\bf Corollary 1.2.} {\it The Fano fiber space $V/{\mathbb P}^1$
is not birationally rigid. The group of birational self-maps
$\mathop{\rm Bir}V\neq\mathop{\rm Bir}(V/{\mathbb P}^1)$ is
infinite.}

{\bf Remark 1.3.} Since the structure of a Fano fiber space
$\pi\colon V\to{\mathbb P}^1$ is not compatible (that is,
fiber-wise) with respect to the structure $q_V$ of the conic
bundle $V/G$, the latter also cannot be birationally rigid. This
agrees with the fact that $V/G$ does not satisfy the Sarkisov
condition:
$$
4K_G+D=-4H_G,
$$
where $D$ is the discriminant divisor of the fiber space $V/G$,
$H_G$ is the class of a hyperplane section of the hypersurface
$G\subset {\mathbb P}$, $K_G=(-M-2+m)H_G$, $D=4lH_G$.

The group $\mathop{\rm Bir}(V/G)$ of birational self-maps,
preserving the structure of a conic bundle $V/G$, is very large.
Accordingly, on the variety $V$ there are a lot of pencils of
rationally connected varieties: it is sufficient to consider all
pencils of the form
$$
\chi_*|F|,\quad \chi\in\mathop{\rm Bir}(V/G).
$$
However, it is unclear, whether there are birational self-maps on
$V$ that are not compatible with the conic bundle $V/G$.

{\bf Conjecture 1.2.} {\it For general $G$, $W$ on the variety
$V$ there is only one structure of a conic bundle, that is,
$q_V\colon V\to G$, whereas the group of birational self-maps of
the variety $V$ and the group of fiber-wise birational self-maps
of the fiber space $V/G$ coincide:}
$$
\mathop{\rm Bir}V=\mathop{\rm Bir}(V/G).
$$
\vspace{0.3cm}

%%%%%%%%%%%%%%%%%%%%%%%%%%%%%%%%%%%%%%%%%%%%%%%%%%%%%%%%%%%%%%
%%%%%%%%%%%%%%%%%%%%%% subsection 1.5

{\bf 1.5. Varieties of type ((0),(0,0)).} Obviously, any variety
of this type is isomorphic to a direct product $F\times {\mathbb
P}^1$, where $\sigma\colon F\to G$ is a Fano double cover of
index 1, branched over the divisor $W_G=G\cap W^*$, $W^*\subset
{\mathbb P}$ is a hypersurface of degree $2l$. Thus the variety
$V$ is not birationally rigid. The morphism of projection onto
the direct factor $q\colon V\to F$ defines on $V$ a conic bundle
structure. Acting, as above, by fiber-wise birational self-maps
$\chi\in \mathop{\rm Bir}(V/F)$ on the pencil $|F|$, we obtain
infinitely many pencils of rationally connected varieties on $V$.
Let $\tau\in \mathop{\rm Aut}V$ be the Galois involution of the
double cover $\sigma\colon V\to Q$, so that for general $Q$, $W$
we have $\mathop{\rm Aut}F={\mathbb Z}/2{\mathbb Z}=\{\mathop{\rm
id},\tau\}$.

{\bf Conjecture 1.3.} {\it For general $G$, $W_*$ there is a
unique conic bundle structure on $V$, that is, $q\colon V\to F$.
For the group of birational self-maps of the variety $V$ the
following presentation holds:}
$$
\mathop{\rm Bir}V=\mathop{\rm Aut}F\times\mathop{\rm
Bir}(V/F)={\mathbb Z}/2{\mathbb Z}\times \mathop{\rm Bir}(V/F).
$$
\vspace{0.3cm}

%%%%%%%%%%%%%%%%%%%%%%%%%%%%%%%%%%%%%%%%%%%%%%%%%%%%%%%%%%%%%
%%%%%%%%%%%%%%%% subsection 1.6

{\bf 1.6. Varieties of type ((1),(0,0)).} First of all, recall
the following well known fact: for $a_1=\dots=a_M=0$, $a_{M+1}=1$
the variety $X={\mathbb P}({\cal E})$ is isomorphic to the blow
up of ${\mathbb P}^{M+2}$ at a subspace of codimension 2, where
the pencil of fibers of the morphism $\pi_X$ corresponds with
respect to this isomorphism to the pencil of hyperplanes in
${\mathbb P}^{M+2}$ containing the centre of the blow up. More
precisely, let $P\subset {\mathbb P}^{M+2}$ be a linear subspace
of codimension two, $\varphi\colon \widetilde {\mathbb P}\to
{\mathbb P}^{M+2}$ its blow up, $E\subset \widetilde {\mathbb P}$
the exceptional divisor, $E\cong {\mathbb P}^M\times{\mathbb
P}^{1}$, whereas $\varphi\colon E\to P={\mathbb P}^{M}$ is just
the projection of $E$ onto the first direct factor. Set ${\cal
L}_P=\varphi^*{\cal O}_{\mathbb P}(1)$. In accordance with our
notations,
$$
{\cal E}={\cal O}^{\oplus (M+1)}_{{\mathbb P}^1}\oplus {\cal
O}_{\mathbb P^1}(1),
$$
where ${\mathbb P}({\cal E})$ is the projectivization of this
sheaf, ${\cal L}$ is the corresponding tautological sheaf,
$E^+\subset{\mathbb P}({\cal E})$ is the divisor of common zeros
of all sections
$$
s\in H^0({\mathbb P}^1,{\cal O}_{{\mathbb P}^1}(1))\hookrightarrow
H^0({\mathbb P}({\cal E}),{\cal L})
$$
(in the sense of the inclusion ${\cal O}_{{\mathbb
P}^1}(1)\hookrightarrow{\cal E}$). Obviously, $E^+={\mathbb
P}({\cal E}^+)$, where ${\cal E}^+={\cal O}^{\oplus
(M+1)}_{{\mathbb P}^1}$, that is, $E^+\cong{\mathbb P}^M\times
{\mathbb P}^{1}$.

{\bf Lemma 1.8.} {\it The varieties $\widetilde{\mathbb P}$ and
${\mathbb P}({\cal E})$ are isomorphic. Moreover, the following
diagram commutes:}
$$
\begin{array}{rcccl}
& \widetilde {\mathbb P} & \leftrightarrow & {\mathbb P}({\cal E}) & \\
\varphi & \downarrow & & \downarrow & \varphi_{\cal L} \\
& {\mathbb P}^{M+2} & = & {\mathbb P}^{M+2}, &
\end{array}
$$
{\it where $\varphi_{\cal L}$ is the morphism determined by the
space of global sections of the tautological sheaf ${\cal L}$.}

{\bf Proof.} This is obvious.

{\bf Lemma 1.9.} {\it For any integer $e\geq 1$ and any
irreducible divisor $R\subset{\mathbb P}({\cal E})$, $R\neq E^+$,
$R=(s)$, where $s\in H^0({\mathbb P}({\cal E}),{\cal L}^{\otimes
e})$, the image $\varphi_{\cal L}(R)\subset {\mathbb P}$ is an
irreducible hypersurface of degree $e$.}

{\bf Proof.} By the arguments above,
$$
{\mathbb P}({\cal E})\setminus E^+\cong {\mathbb
P}^{M+2}\setminus P,
$$
where the tautological divisors on ${\mathbb P}({\cal E})$
correspond to hyperplanes on ${\mathbb P}^{M+2}$. A general line
$L\subset{\mathbb P}^{M+2}$ does not meet $P$ and thus
$$
\mathop{\rm deg}({\cal L}\,|_{\varphi^{-1}_{\cal L}(L)})=1.
$$
In other words, $(R\cdot\varphi^{-1}_{\cal L}(L))=(\varphi_{\cal
L}(R)\cdot L)=e$, which is what we need.

By Lemma 1.9, the variety $V$ is birational to the variety
$V^{\sharp}$ which is realized as the double cover
$\sigma^{\sharp}\colon V^{\sharp}\to Q^{\sharp}$,
$Q^{\sharp}=\varphi_{\cal L}(Q)\subset {\mathbb P}^{M+2}$, and
$\sigma^{\sharp}$ is branched (in codimension 1) over the divisor
$Q^{\sharp}\cap W^{\sharp}$, where $W^{\sharp}=\varphi_{\cal
L}(W)\subset{\mathbb P}^{M+2}$. By Lemma 1.9, $\mathop{\rm
deg}Q^{\sharp}=m$, $\mathop{\rm deg}W^{\sharp}=2l$, so that
$V^{\sharp}$ is a Fano double hypersurface of index 2. Any pencil
of hyperplanes in ${\mathbb P}^{M+2}$ determines a pencil of Fano
varieties on $V^{\sharp}$ and thus a pencil of rationally
connected divisors on $V$.

{\bf Conjecture 1.4.} {\it Apart from the pencils of rationally
connected divisors described above, there are no other structures
of a rationally connected fiber space on $V$. The groups of
birational and biregular automorphisms of the variety $V$
coincide:}
$$
\mathop{\rm Bir}V=\mathop{\rm Aut}V={\mathbb Z}/2{\mathbb Z}.
$$
\vspace{0.3cm}

%%%%%%%%%%%%%%%%%%%%%%%%%%%%%%%%%%%%%%%%%%%%%%%%%%%%%%%%%%%%%%%%
%%%%%%%%%%%%%%%%%% subsection 1.7

{\bf 1.7. Varieties with a pencil of double spaces.} Recall that
when Fano double spaces of index 1 are considered, the symbol
${\mathbb P}$ denotes the projective space ${\mathbb P}^M$ of one
dimension less. \vspace{0.3cm}

{\bf Varieties of type ((0),1).} Here $X={\mathbb
P}\times{\mathbb P}^1$ and the branch divisor $W\subset X$ is of
bidegree $(2M,2)$, that is, it is given by the equation
$$
f_{uu}u^2+2f_{uv}uv+f_{vv}v^2=0,
$$
where $(u:v)$ are homogeneous coordinates on ${\mathbb P}^1$,
$f_{\sharp}$ are homogeneous polynomials of degree $2M$ on
${\mathbb P}$. This case is completely similar to the one studied
above in Sec. 1.4. The projection $q\colon V\to {\mathbb P}$
realizes $V$ as a conic bundle. Its discriminant divisor
$D\subset{\mathbb P}$ is a hypersurface of degree $4M$, that is,
the Sarkisov condition is not satisfied.

{\bf Conjecture 1.5.} {\it On the variety $V$ there are no other
structures of a conic bundle, apart from $V/{\mathbb P}$. The
groups of birational and fiber-wise birational self-maps coincide:
}
$$
\mathop{\rm Bir}V=\mathop{\rm Bir}(V/{\mathbb P}).
$$

{\bf Varieties of type ((0),0).} Here $X={\mathbb
P}\times{\mathbb P}^1$, $W=W_{\mathbb P}\times {\mathbb P}^1$,
where $W_{\mathbb P}\subset {\mathbb P}$ is a hypersurface of
degree $2M$. Let $\sigma\colon F\to{\mathbb P}$ be the double
space branched over the hypersurface $W_{\mathbb P}$. This case
is completely similar to the one considered above in Sec. 1.5:
$V\cong F\times {\mathbb P}^1$.

{\bf Conjecture 1.6.} {\it Apart from $V/F$, there are no other
structures of a conic bundle on $V$. For the group of birational
self-maps the following presentation holds:}
$$
\mathop{\rm Bir}V=\mathop{\rm Aut}F\times \mathop{\rm
Bir}(V/F)={\mathbb Z}/2{\mathbb Z}\times \mathop{\rm Bir}(V/F).
$$
\vspace{0.3cm}

{\bf Varieties of type ((1),0).} This case is completely similar
to the case considered above in Sec. 1.6. Here the variety
$X={\mathbb P}({\cal E})$ is isomorphic to the blow up of
${\mathbb P}^{M+1}$ at a linear space $P\subset{\mathbb P}^{M+1}$
of codimension two. Elementary computations (see Sec. 1.6) show
that with respect to this isomorphism the variety $V$ is
birationally equivalent to a double space ${\mathbb P}^{M+1}$ of
index 2 (with the branch divisor $W_{\mathbb P}\subset {\mathbb
P}^{M+1}$ of degree $2M$). Moreover, it is easy to see that all
varieties of type ((1),0) are realized in this way: take a double
space of index 2 and any pencil of hyperplanes on ${\mathbb
P}^{M+1}$. In particular, to any plane of codimension two
$P\subset{\mathbb P}^{M+1}$ corresponds a pencil $\Lambda_P$ of
Fano varieties on $V$.

{\bf Conjecture 1.7.} {\it Every structure of a rationally
connected Fano fiber space on $V$ is a pencil $\Lambda_P$ for
some subspace $P\subset {\mathbb P}^{M+1}$of codimension two. The
groups of birational and biregular automorphisms of the variety
$V$ coincide:}
$$
\mathop{\rm Bir}V=\mathop{\rm Aut}V.
$$
\vspace{0.3cm}

{\bf Varieties of type ((1,1),0).} This case is completely
similar to the case considered above in Sec. 1.2. For these
varieties the claim of Proposition 1.1 holds. The proof given in
Sec. 1.2 works with some simplifications (there is no divisor
$Q$).

%%%%%%%%%%%%%%%%%%%%%%%%%%%%%%%%%%%%%%%%%%%%%%%%%%%%%%%%%%%%%%%%%
%%%%%%%%%%%%%%%%%%%%%%%%%%%%%%%%%%%%%%%%%%%%%%%%%%%%%%%%%%%%%%%%%
%%%%%%%%%%%%%%%% SECTION TWO

\section{Infinitely near singularities of vertical subvarieties}

{\bf 2.1. Set up of the problem.} Recall that in [1, Sec. 2.2]
the condition (vs) for the case of a singular point of a fiber
$o\in F$ lying outside the ramification divisor of the double
cover $\sigma\colon F\to G$ was proved in the following way.
Assume that there exists a prime divisor $Y\subset F$, satisfying
the estimate
\begin{equation}
\label{b1} \frac{\mathop{\rm mult}_x\widetilde Y}{\mathop{\rm
deg}Y}>\frac{1}{m},
\end{equation}
where $x$ is an infinitely near point of the first order over the
point $o$, that is, $x\in E$, where $E\subset\widetilde F$ is the
exceptional divisor of the blow up $\varphi\colon \widetilde F\to
F$ of the point $o$. Here $p=\sigma(o)\in G$ is a non-degenerate
double point of the hypersurface $G$, $p\not\in W$.

Now the crucial fact is Proposition 2.2 in [1, Sec. 2.2]:

{\it There exists a hyperplane $P\subset {\mathbb P}$, $P\ni p$,
such that $\sigma(Y)\not\subset P$ and the effective algebraic
cycle $Y_P=(Y\mathop{\circ}\nolimits_F P_F)$, where
$P_F=\sigma^{-1}(P_G)$, $P_G=P\cap G$ is the corresponding
hyperplane section, satisfies the estimate
\begin{equation}\label{b2}
\frac{\mathop{\rm mult}\nolimits_o}{\mathop{\rm deg}} Y_P
>\frac{3}{2m}.
\end{equation}
}

{\bf Proof} repeats the arguments in the similar case in [7] word
for word and for this reason was omitted in [1]. It is based on
the following local fact. Set $2\nu=\mathop{\rm mult}\nolimits_o
Y$, $\mu=\mathop{\rm mult}\nolimits_x\widetilde Y$, $B=T_xE\cap
E$, where $E$ is considered as a quadratic hypersurface in
${\mathbb P}^M$.

{\bf Lemma 2.1.} {\it The following estimate holds}
\begin{equation}\label{b1a}
\mathop{\rm mult}\nolimits_B\widetilde Y\geq \frac12(\mu-\nu).
\end{equation}
Proposition 2.2 of the paper [1] cited above follows from Lemma
2.1 almost immediately. Indeed, the exceptional divisor $E$ is
embedded in the exceptional divisor ${\mathbb T}$ of the blow up
of the point $p$ on the projective space ${\mathbb P}$ as a
quadric hypersurface. Thus one may consider the tangent
hyperplane $T_xE$ as a hyperplane in ${\mathbb T}$, that is, as
the projectivized tangent cone to a uniquely determined hyperplane
$P\subset{\mathbb P}$, $p\in P$. Let $P^{\sharp}$ be the strict
transform of $P$ on $\widetilde{\mathbb P}$, then
$$
P^{\sharp}\cap{\mathbb T}=T_xE.
$$
Set $P_G=P\cap G$ and $P_F=\sigma^{-1}(P_G)$. These are
irreducible varieties and moreover
$$
\mathop{\rm deg}P_F=2m, \quad \mathop{\rm
mult}\nolimits_x\widetilde P_F\leq 2,
$$
so that $Y\neq P_F$. Furthermore,
$$
B=T_xE\cap E\subset\widetilde P_F,
$$
so that by the standard formulas of intersection theory [8,
Chapter 2] we get for the effective cycle
$Y_P=(Y\mathop{\circ}_FP_F)$:
$$
\mathop{\rm mult}\nolimits_oY_P\geq 2\nu+(\mathop{\rm
mult}\nolimits_B\widetilde Y)\mathop{\rm deg}B= \mu+\nu,
$$
whence, taking into account that $\mathop{\rm deg}
Y_P=\mathop{\rm deg}Y$ (since $Y_P$ is a hyperplane section of
$Y$) and by assumption (see (\ref{b1}))
$$
2\nu\geq \mu=\mathop{\rm mult}\nolimits_x \widetilde
Y>\frac{1}{m}\mathop{\rm deg}Y,
$$
we obtain the desired inequality (\ref{b2}).

Unfortunately, this simple argument, by means of which
Proposition 2.2 of the paper [1] is derived from Lemma 2.1, does
not work in the case when the singular point $o\in F$ lies on the
ramification divisor: if, in the notations of Sec. 3.1 of [1],
the point $x\in E$ lies outside the ramification divisor of the
double cover $\tilde\sigma_E$, that is, $\tilde\sigma(x)\not\in
W_E$, then the hyperplane $P\subset{\mathbb P}$, described above,
does not exist. The best we could take instead of the subvariety
$P$ in this case is a quadric hypersurface, the tangent cone to
which at the point $p$ contains the quadric $\tilde\sigma_E(B)$.
However, the ratio of the multiplicity to the degree for the
cycle $Y_P$ turns out to be not big enough. That is why in [1] an
alternative method of proving the condition (vs) was developed
for this case. Let us consider in a more detailed way, how the
lower estimate for the number $(\mathop{\rm
mult}\nolimits_o/\mathop{\rm deg})Y_P$ is obtained. Let
$P\subset{\mathbb P}$ be a quadric hypersurface, containing the
point $p\in W\cap G$, where $p\in P$ is a singular point and
$P^{\sharp}\supset\tilde\sigma_E(B)$. Now repeating the arguments
above word for word, we get the estimate
\begin{equation}\label{b3}
\mathop{\rm mult}\nolimits_o Y_P\geq4\nu+2\mathop{\rm
mult}\nolimits_B \widetilde Y,
\end{equation}
whence, taking into account Lemma 2.1 we get
$$
\frac{\mathop{\rm mult}\nolimits_o}{\mathop{\rm
deg}}Y_P>\frac{5}{4m},
$$
but this is not good enough to get a contradiction. However, it
is easy to see from the inequality (\ref{b3}) that what we need
to obtain an estimate for the multiplicity of the cycle $Y_P$ at
the point $o$, is not so much an estimate of the number
$\mathop{\rm mult}\nolimits_B \widetilde Y$, but rather a
combined estimate for $\nu$ and $\mathop{\rm
mult}\nolimits_B\widetilde Y$. For instance, if we knew that
\begin{equation}
\label{b4} \nu+\mathop{\rm mult}\nolimits_B \widetilde
Y>\frac{1}{m}\mathop{\rm deg}Y,
\end{equation}
then, taking into account that $\nu\geq\mathop{\rm
mult}\nolimits_B \widetilde Y$, we could get
$$
\frac{\mathop{\rm mult}\nolimits_o}{\mathop{\rm
deg}}Y_P>\frac{3}{2m},
$$
after which we could have argued as in the case of a singular
point outside the ramification divisor. The author believes that
if (in the notations of [2, Sec. 1.4]) the inequality $k_v>2$
holds, then the estimate (\ref{b4}) is true. The aim of this
section is to formulate the corresponding claim precisely and to
prove it in the particular case when the graph of the sequence of
blow ups is a chain. Proof of this conjecture in full would have
essentially simplified our work with Fano fiber spaces
$V/{\mathbb P}^1$ in the case of vertical subvarieties.

The structure of this section is as follows. To begin with, we
recall the proof of Lemma 2.1 (following [7], but with much more
details). We need it to develop new arguments generalizing the
method of proof of Lemma 2.1. After that we consider the general
case and formulate the above-mentioned conjecture. In the
remaining part of the section we prove this conjecture for the
case when the graph of the sequence of blow ups is a chain.
\vspace{0.3cm}

%%%%%%%%%%%%%%%%%%%%%%%%%%%%%%%%%%%%%%%%%%%%%%%%%%%%%%%%%%%%%%%%%%%
%%%%%%%%%%%%%%%%%%%%%%%%%%%%%%%%%%%%%%%%%%%%%%%%%%%%%%%%%%%%%%%%%%%
%%%%%%%%%%%%%%%%%%%%%%%%% subsection 2.2  %%%%%%%%%%%%%%%%%%%%%%%%%

{\bf 2.2. Proof of Lemma 2.1.} The claim of the lemma is local.
Let $\Pi\ni p$ be a germ of a general section of the fiber $G$ by
a smooth four-dimensional subvariety of the ambient space
(analytically, $G$ in a neighborhood of the point $p$ is a
hypersurface in ${\mathbb C}^{M+1}$), such that its strict
transform on $\widetilde G$ contains the point $x$: $\widetilde
\Pi\ni x$. Set
$$
E\cap\widetilde\Pi=E_{\Pi}\cong {\mathbb P}^1\times {\mathbb P}^1.
$$
Obviously,
$$
B_{\Pi}=B\cap E_{\Pi}=L_1+L_2
$$
is the union of the two lines on the quadric $E$, that pass
through the point $x$. Since $\Pi\ni p$ is a germ of
three-dimensional non-degenerate quadratic singularity, its
strict transform $\widetilde\Pi$ is smooth. Set $L=L_1$.

Let
$$
Y_{\Pi}=Y\cap \Pi, \quad \widetilde Y_{\Pi}=\widetilde Y\cap
\widetilde \Pi
$$
be the restrictions of $Y$ onto $\Pi$, $\widetilde \Pi$
respectively. Let us prove the inequality
$$
\gamma=\mathop{\rm mult}\nolimits_L\widetilde Y_{\Pi}\geq
\frac12(\mu-\nu).
$$
Since the germ $\Pi$ is a general one, this implies the inequality
(\ref{b1a}).

Let
$$
\varphi_L\colon \Pi_L\to\widetilde\Pi
$$
be the blow up of the line $L$ on the smooth three-dimensional
variety $\widetilde \Pi$, $E_L\subset \Pi_L$ the exceptional
divisor. Since
$$
(L\cdot L)_{E_{\Pi}}=1, \quad (L\cdot E_{\Pi})_{\widetilde
\Pi}=-1,
$$
we obtain the exact sequence
$$
\begin{array}{rcccl}
0\to & {\cal N}_{L/E_{\Pi}} & \to {\cal N}_{L/\widetilde \Pi}\to &
{\cal N}_{E_{\Pi}/{\widetilde \Pi}}|_L & \to 0, \\
& \| & & \|&  \\
& {\cal O}_L & & {\cal O}_L(-1)&
\end{array}
$$
so that
\begin{equation}\label{b5}
{\cal N}_{L/\widetilde \Pi}\cong{\cal O}_L\oplus{\cal O}_L(-1)
\end{equation}
and thus $E_L$ is a ruled surface of type ${\mathbb F}_1$ with
the exceptional section $\widetilde E_{\Pi}\cap E_L$, the class
of which is denoted by $s$. We have
$$
\mathop{\rm Pic}E_L={\mathbb Z}s\oplus{\mathbb Z}f,
$$
where $f$ is the class of a fiber of the ruled surface. Since
$$
(s\cdot E_L|_{E_L})_{E_L}=((\widetilde E_{\Pi}\cap E_L)\cdot
 E_L)=(L\cdot L)_{L_{\Pi}}=0,
$$
we get $(E_L\cdot E_L)=-s-f$ (which can also be seen directly from
(\ref{b5})). Let $Y_L$ be the strict transform of the divisor
$\widetilde Y_{\Pi}$ on $\Pi_L$. Obviously,
$$
Y_L\sim-\nu E_{\Pi}-\gamma E_L
$$
(recall that $\Pi$ is a germ of a three-dimensional section).
Therefore,
$$
Y_L|_{E_L}\sim \gamma s+(\gamma+\nu)f.
$$
On the other hand, the following fact is well known.

{\bf Lemma 2.2.} {\it Let $Z\subset R\subset X$ be a flag of
strictly embedded smooth varieties, $D$ an effective divisor on
$X$,
$$
\varphi_R\colon \widetilde X\to X
$$
the blow up of the subvariety $R$ with the exceptional divisor
$E_R$ and $E_Z=\varphi^{-1}_R(Z)$. Let $\widetilde D\subset
\widetilde X$ be the strict transform of the divisor $D$ on
$\widetilde X$. The following estimate holds:}
$$
\mathop{\rm mult}\nolimits_{E_Z}\widetilde D=\mathop{\rm
mult}\nolimits_Z D-\mathop{\rm mult}\nolimits_RD.
$$

{\bf Proof.} Restricting $Z$, $R$, $D$ onto a general smooth
subvariety in $X$ of a suitable dimension, we assume that
$Z=\{z\}$ is a point, $R$ is a curve. Let $S\supset R$ be a
general surface, $\widetilde S\cong S$ its strict transform on
$\widetilde X$. We get
$$
D|_S=(\mathop{\rm mult}\nolimits_R D)R+D^{\sharp}_S,
$$
$$
\widetilde D|_{\widetilde S}=D^{\sharp}_S,\quad \mathop{\rm
mult}\nolimits_{E_Z}\widetilde D=\mathop{\rm mult}\nolimits_z
D^{\sharp}_S,
$$
but since $\mathop{\rm mult}_Z D=\mathop{\rm mult}_z D|_S$ and
$\mathop{\rm mult}_z R=1$, we obtain the claim of the lemma.

Let $\Lambda\subset E_L$ be the fiber over the point $x$. By the
lemma just proved, we get the estimate
$$
\mathop{\rm mult}\nolimits_{\Lambda}(Y_L|_{E_L})\geq \mu-\gamma.
$$
Therefore,
$$
\gamma+\nu\geq \mu-\gamma,
$$
so that $2\gamma\geq \mu-\nu$, which is what we need.
\vspace{0.3cm}

%%%%%%%%%%%%%%%%%%%%%%%%%%%%%%%%%%%%%%%%%%%%%%%%%%%%%%%%%%%%%%%%%%%
%%%%%%%%%%%%%%%%%%%%%%%%%%%%%%%%%%%%%%%%%%%%%%%%%%%%%%%%%%%%%%%%%%%
%%%%%%%%%%%%%%%%%%%%%%%%% subsection 2.3  %%%%%%%%%%%%%%%%%%%%%%%%%

{\bf 2.3. A conjecture on multiplicities.} Let $\Delta\ni o$ be a
germ of a non-degenerate quadratic singularity, $\mathop{\rm dim}
\Delta=N+1$, $N\geq 3$. Let us blow up the point $o$:
$$
\varphi\colon \widetilde \Delta\to\Delta,
$$
$E\subset\widetilde\Delta$ is the exceptional divisor, that is, an
$N$-dimensional quadric. Let
$$
\varphi_{i,i-1}\colon \Delta_i\to\Delta_{i-1}
$$
be a sequence of blow ups of irreducible varieties
$B_{i-1}\subset \Delta_{i-1}$, $i=1,\dots,k$, $\Delta_0=\widetilde
\Delta$, with the exceptional divisors
$$
E_i=\varphi^{-1}_{i,i-1}(B_{i-1})\subset\Delta_i,
$$
where, moreover, $B_i\subset E_i$ for $i=0,\dots,k-1$, and
$E_0=E$, $B_0=x\in E$ is a point. Consider the graph $\Gamma$,
associated with the sequence of blow ups $\varphi_{i,i-1}$, and
set
$$
p_i=\sharp\{ \mbox{the paths from}\,\, E_k \,\,\mbox{to}\,\, E_i
\}.
$$
Let $D\subset\Delta$ be a prime divisor,
$$
D^0=\widetilde D\subset \widetilde\Delta=\Delta_0,\quad
D^1\subset\Delta_1,\dots,D^k\subset \Delta_k
$$
its strict transforms on $\Delta_0=\widetilde\Delta$,
$\Delta_1,\dots,\Delta_k$, respectively. We get
$$
D^k\sim D-\nu E-\mu_1 E_1-\dots-\mu_k E_k,
$$
for some  $\nu,\mu_1,\dots,\mu_k\in{\mathbb Z}_+$. More precisely,
$$
2\nu=\mathop{\rm mult}\nolimits_o D,\quad \mu_i=\mathop{\rm
mult}\nolimits_{B_{i-1}}D^{i-1}.
$$
Set
$$
B=T_xE\cap E,
$$
where the quadric $E$ is assumed to be embedded in ${\mathbb
P}^{N+1}$ in the standard way.

{\bf Conjecture 2.1.} {\it Assume that the following estimate
holds:
$$
2\nu p_0+\sum^k_{i=1}p_i\mu_i\geq \lambda \left(
2p_0+\sum^k_{i=1}p_i \right).
$$
Then the following inequality is satisfied:}
\begin{equation}
\label{b6} \nu+\mathop{\rm mult}\nolimits_B \widetilde D\geq
\lambda.
\end{equation}
We prove this conjecture for the particular case when the graph
$\Gamma$ is a chain, that is,
$$
B_i\not\subset E^i_{i-1}
$$
for every $i=2,\dots,k-1$. In this case all the integers $p_i$
are equal to 1.\vspace{0.3cm}

%%%%%%%%%%%%%%%%%%%%%%%%%%%%%%%%%%%%%%%%%%%%%%%%%%%%%%%%%%%%%%%%%%%
%%%%%%%%%%%%%%%%%%%%%%%%%%%%%%%%%%%%%%%%%%%%%%%%%%%%%%%%%%%%%%%%%%%
%%%%%%%%%%%%%%%%%%    subsection 2.4   %%%%%%%%%%%%%%%%%%%%%%%%%%%%
{\bf 2.4. Proof in the case when the graph is a chain.}

{\bf Proposition 2.1.} {\it Assume that the graph $\Gamma$ is a
chain and the inequality
\begin{equation}
\label{b7} 2\nu+\mu_1+\dots+\mu_k\geq (k+2)\lambda
\end{equation}
holds. Then the estimate (\ref{b6}) is true.}

{\bf Proof} is given in the assumption that all the centres of
blow ups $B_i$ are points. The general case reduces to this one
in a trivial way. We get
$$
E_1\cong \dots\cong E_k\cong {\mathbb P}^N.
$$
Let us construct by induction a sequence of irreducible
non-singular subvarieties of codimension two
$$
Y_j\subset \Delta_{j+2},
$$
$j=0,\dots,k-1$, in the following way. Set $Y_0=B^1$ to be the
strict transform of the quadric cone $B$ on $\Delta_1$. Obviously,
$Y_0$ is the cone $B$ with the blown up vertex; in particular,
$Y_0\cap E_1$ is a non-singular quadratic hypersurface in the
hyperplane $E^1\cap E_1\subset E_1\cong {\mathbb P}^N$. Now
$$
Y_1,\dots,Y_{k-1}
$$
are uniquely determined by the following condition: the subvariety
$$
Y^+_j=\varphi_{j+1,j}(Y_j)\subset E_j
$$
is a non-degenerate quadratic cone in $E_j\cong {\mathbb P}^N$
with the vertex at the point $x_j\in E_j$ and the base
$Y_{j-1}\cap E_j$, and moreover $Y_j$ is the strict transform of
$Y^+_j$ on $\Delta_{j+1}$.

This construction is justified exactly by the fact that the
non-singular quadric $Y_{j-1}\cap E_j$ is contained in the
hyperplane $E^j_{j-1}\cap E_j$, whereas the vertex $x_j$ of the
cone is not contained in this hyperplane by assumption:
$x_j\not\in E^j_{j-1}$, since the graph $\Gamma$ is a chain.

Set
$$
\gamma_i=\mathop{\rm mult}\nolimits_{Y_i} D^i=\mathop{\rm
mult}\nolimits_{Y^+_i}D^{i-1}.
$$
For convenience set also $\gamma_{-1}=0$, $\mu_0=\nu$,
$\gamma_k=0$.

{\bf Lemma 2.3.} {\it For any $i\in\{0,\dots,k-1\}$ the following
estimate holds:}
\begin{equation}
\label{b8} 2\gamma_i\geq
\mu_{i+1}-\mu_i+\gamma_{i-1}+\gamma_{i+1}.
\end{equation}

{\bf Proof} of the lemma is given below. Now let us complete the
proof of Proposition 2.1. The estimates (\ref{b8}), put together
with the coefficients $(k-i)$, give the inequality
\begin{equation}
\label{b9}
\sum^{k-1}_{i=0}(k-i)(2\gamma_i-\mu_{i+1}+\mu_i-\gamma_{i-1}-\gamma_{i+1})\geq
0.
\end{equation}
For $i\neq 0$ the multiplicity $\gamma_i$ comes into the
left-hand side of (\ref{b9}) with the coefficient
$$
2(k-i)-(k-i-1)-(k-i+1)=0;
$$
the multiplicity $\gamma_0$ comes into the left-hand side of
(\ref{b9}) with the coefficient
$$
-(k-1)+2k=k+1;
$$
for $i\neq 0$ the multiplicity $\mu_i$ comes with the coefficient
$$
(k-i)-(k-i+1)=-1;
$$
finally, the multiplicity $\mu_0$ comes into (\ref{b9}) with the
coefficient $k$. As a result, we get the inequality
$$
(k+1)\gamma_0+k\nu\geq \mu_1+\dots+\mu_k,
$$
whence, taking into account the inequality (\ref{b7}), we get
$$
(k+1)\gamma_0+(k+2)\nu\geq (k+2)\lambda,
$$
so that the more so,
$$
\gamma_0+\nu\geq\lambda.
$$
Now recall that $\gamma_0=\mathop{\rm mult}\nolimits_B \widetilde
D$, which completes the proof of Proposition 2.1. \vspace{0.3cm}

%%%%%%%%%%%%%%%%%%%%%%%%%%%%%%%%%%%%%%%%%%%%%%%%%%%%%%%%%%%%%%%%%%%
%%%%%%%%%%%%%%%%%%%%%%%%%%%%%%%%%%%%%%%%%%%%%%%%%%%%%%%%%%%%%%%%%%%
%%%%%%%%%%%%%%%%%%    subsection 2.5   %%%%%%%%%%%%%%%%%%%%%%%%%%%%
{\bf 2.5. Proof of Lemma 2.3.} Let us consider the following
general situation. Let $\Pi\ni o$ be a smooth three-dimensional
germ,
$$
\psi\colon \widetilde \Pi\to\Pi
$$
the blow up of the point $o$, $E\subset\widetilde \Pi$ the
exceptional divisor, $E\cong{\mathbb P}^2$. Let $a\in E$ be an
arbitrary point,
$$
\psi^+\colon \Pi^+\to \widetilde\Pi
$$
its blow up, $E^+\subset \Pi^+$ the exceptional plane, $\widetilde
E$ the strict transform of the plane $E$ on $\Pi^+$. Obviously,
$\widetilde E$ is a ruled surface of type ${\mathbb F}_1$.
Finally, let $L^*\subset E$ be a line, passing through the point
$a$, $L\subset \widetilde E$ its strict transform on $\Pi^+$, and
let $R\subset E^+$ be a line, passing through the point $ L\cap
E^+$, but different from the line $\widetilde E\cap E^+$. On the
line $L$ we take a point $b\neq a$, and let $C\subset \Pi$ be a
smooth curve, intersecting transversally the plane $E$ at the
point $b$.

Let $D\subset \Pi$ be a germ of a prime divisor, $o\in D$. Set
$\widetilde D\subset \widetilde\Pi$ and $D^+\subset \Pi^+$ to be
the strict transforms of $D$ on $\widetilde \Pi$ and $\Pi^+$,
respectively. Finally, set
$$
\mu_o=\mathop{\rm mult}\nolimits_o D,\quad \mu_a=\mathop{\rm
mult}\nolimits_a \widetilde D,
$$
$$
\mu_C=\mathop{\rm mult}\nolimits_C \widetilde D,\quad
\mu_R=\mathop{\rm mult}\nolimits_R D^+,
$$
$$
\mu_L=\mathop{\rm mult}\nolimits_L \widetilde D.
$$
The following claim holds.

{\bf Lemma 2.4.} {\it We have the estimate}
$$
2\mu_L\geq \mu_a+\mu_R+\mu_C-\mu_o.
$$

{\bf Proof.} Blow up the curve $L$:
$$
\psi_L\colon \Pi_L\to \Pi^+,
$$
and let $E_L\subset \Pi_L$ be the exceptional divisor. It is well
known that the normal sheaf of the line $L^*\subset \widetilde
\Pi$ is of the form
$$
{\cal N}_{L^*/\widetilde \Pi}\cong {\cal O}_{L^*}(-1)\oplus {\cal
O}_{L^*}(1),
$$
so that the normal sheaf of the curve $L$ is
\begin{equation}
\label{b10} {\cal N}_{L/\Pi^+}\cong {\cal O}_L(-2)\oplus {\cal
O}_L.
\end{equation}
In particular, $E_L$ is a surface of type ${\mathbb F}_2$. Set
$E^L$ to be the strict transform of the surface $\widetilde E$ on
$\Pi_L$. It is easy to see that the exceptional section of the
ruled surface $E_L$ is the curve $E_L\cap E^L$. Write down
$$
\mathop{\rm Pic} E_L={\mathbb Z} s\oplus {\mathbb Z} f,
$$
where $(s^2)=-2$, $f$ is the class of a fiber. Denote by the
symbols $C_L$ and $R_L$ the strict transforms of the curves $C$
and $R$ on $\Pi_L$. Since $C$ and $R$ intersect the surface
$\widetilde E$ transversally at the points $b\in L$ and $L\cap
E^+$, the curves $C_L$ and $R_L$ intersect (transversally) the
ruled surface $E_L$ at the points $x_C=C_L\cap E_L$ and
$x_R=R_L\cap E_L$, which do not lie on the exceptional section
$E_L\cap E^L$, respectively. Finally, let $D_L\subset \Pi_L$ be
the strict transform of the divisor $D$. We get
$$
D_E=(D_L\circ E_L)\sim \mu_L (-E_L|_{E_L})+(\mu_o-\mu_a)f,
$$
since, obviously,
$$
(\psi_L\circ \psi^+)^* E|_{E_L}\sim -f, \quad \psi^*_L
E^+|_{E_L}\sim f.
$$
From the form of the normal sheaf (\ref{b10}) or from the relation
$$
(E-E^+-E_L)|_{E_L}\sim s
$$
(since  $E^L\sim E-E^+-E_L$) we obtain that
$$
-E_L|_{E_L}\sim s+2f,
$$
so that we get finally
\begin{equation}
\label{b11} D_E\sim \mu_L s+(2\mu_L+\mu_o-\mu_a)f.
\end{equation}
However, the curve  $D_E\subset E_L$ is effective, and moreover,
\begin{equation}
\label{b12} \mathop{\rm mult}\nolimits_{x_C}D_E\geq \mu_C,\quad
\mathop{\rm mult}\nolimits_{x_R}D_E\geq \mu_R.
\end{equation}
Write down
\begin{equation}
\label{b13} D_E=A^*+\alpha_C A_C+\alpha_R A_R,
\end{equation}
where $A_C$ (respectively, $A_R$) is the fiber of the ruled
surface $E_L$, containing the point $x_C$ (respectively, $x_R$),
and the curve $A^*$ is effective and does not contain $A_C$ or
$A_R$ as a component. Obviously,
$$
A^*\sim \mu_L s+\alpha^* f,
$$
where $\alpha^*\in{\mathbb Z}_+$.

{\bf Lemma 2.5.} {\it The following estimate holds}
\begin{equation}
\label{b14} \mathop{\rm mult}\nolimits_{x_C}A^*+\mathop{\rm
mult}\nolimits_{x_R}A^*\leq \min (\alpha^*,2\mu_L).
\end{equation}

{\bf Proof.} Since $A_C$ is not a component of the curve $A^*$,
we get
\begin{equation}
\label{b15} \mathop{\rm mult}\nolimits_{x_C}A^*\leq (A^*\cdot
A_C)_{x_C}\leq (A^*\cdot A_C)=\mu_L,
\end{equation}
and similarly for $x_R$. Furthermore, if $S$ is an irreducible
component of the curve $A^*$, containing at least one of the
points $x_C,x_R$, then $S$ is not the exceptional section of the
surface $E_L$ (since $x_C,x_R\not\in E_L\cap E^L$). Therefore,
$$
S\sim \beta s+\delta f,
$$
where $\delta\geq 2\beta$, whence as above in (\ref{b15}) we get
$$
\mathop{\rm mult}\nolimits_{x_C} S\leq\beta\leq \frac{\delta}{2}
$$
and similarly for $x_R$, so that
$$
\mathop{\rm mult}\nolimits_{x_C} S+\mathop{\rm
mult}\nolimits_{x_R} S\leq \delta.
$$
This proves Lemma 2.5.

Taking into account all the estimates (\ref{b11} - \ref{b14}) and
the fact that $x_C\not\in A_R$ and $x_R\not\in A_C$ by
construction, we get
$$
2\mu_L+\mu_o-\mu_a\geq \mu_C+\mu_R,
$$
as we claimed. Proof of Lemma 2.4 is complete.

Let us come back to the proof of Lemma 2.3. Let us prove first
the inequality (\ref{b8}) for $i\geq 1$. Let
$$
\Gamma^-_i\subset \Delta_{i-1}
$$
be a germ of a smooth curve, satisfying the following conditions:

(i) $x_{i-1}\in\Gamma^-_i$;

(ii) the strict transform $\Gamma_i\subset\Delta_i$ of the curve
$\Gamma^-_i$ contains the point $x_i$;

(iii) the strict transform $\Gamma^+_i\subset \Delta_{i+1}$ of
the curve $\Gamma_i$ contains the point $x_{i+1}$.

Such a germ exists exactly because the graph
$$
E_i\leftarrow E_{i+1}\leftarrow E_{i+2}
$$
is a chain. Now let $\Pi^-_i\supset \Gamma^-_i$ be a general
three-dimensional germ at the point $x_{i-1}$. Obviously, the
strict transform $\Pi_i\subset \Delta_i$ of the germ $\Pi^-_i$
contains the point $x_i$, whereas, in its turn, its strict
transform $\Pi^+_i$ on $\Delta_{i+1}$ contains the point
$x_{i+1}$. We denote the plane
$$
\Pi_i\cap E_i\subset \Delta_i
$$
by the symbol  $E_{\Pi}$. Since $E_{\Pi}$ contains the vertex of
the quadratic cone $Y^+_i$, the intersection $E_{\Pi}\cap Y^+_i$
is a pair of distinct lines passing through the point $x_i$. Set
$L^*_i\subset E_{\Pi}\cap Y^+_i$ to be one of these two lines,
$L_i\subset \Delta_{i+1}$ its strict transform on $\Delta_{i+1}$.
Let $R_i\subset \Pi^+_i$ be the line in $E_{i+1}\cong {\mathbb
P}^N$, joining the points $x_{i+1}$ and $L_i\cap E_{i+1}$.
Finally, if the germ $\Pi^-_i$ is sufficiently general, then the
curve $Y^+_{i-1}\cap\Pi^-_i$ at the point $x_{i-1}$ has a pair of
branches with distinct tangent directions, so that the curve
$$
C_i=\Pi_i\cap Y_{i-1}\subset E^i_{i-1}
$$
intersects $E_i$ transversally at two distinct points. Let
$a_i\in C_i\cap E_i$ be that one which lies on the line $L^*_i$
(recall that $Y_{i-1}\cap E_i$ is the base of the quadratic cone
$Y^+_i$, and $\Pi_i\cap E_i$ is a plane, containing its vertex).
Now we apply the lemma with $x_{i-1}$, $x_i$, $L_i$, $R_i$, $C_i$
instead of $o$, $a$, $L$, $R$, $C$, respectively. This
immediately implies the estimate (\ref{b8}) for $i\geq 1$.

For $i=0$ the situation is completely similar with the only
exception: there is no curve $C_i$, so that the inequality of
Lemma 2.4 is used in the truncated form:
$$
2\mu_L \geq \mu_a+\mu_R-\mu_o,
$$
where $o\in\Pi$ is a germ of a non-degenerate quadratic
singularity,
$$
\mu_o=\frac12 \mathop{\rm mult}\nolimits_o D,
$$
as in the proof of Lemma 2.1. This completes the proof of Lemma
2.3.

%%%%%%%%%%%%%%%%%%%%%%%%%%%%%%%%%%%%%%%%%%%%%%%%%%%%%%%%%%%%%%%%%%%%
%%%%%%%%%%%%%%%%%%%%%%%%%%%%%%%%%%%%%%%%%%%%%%%%%%%%%%%%%%%%%%%%%%%%
%%%%%%%%%%%%%%%%%%%%  SECTION 3: LINEAR METHOD

\section{The linear method}

{\bf 3.1. The linear method of proving birational rigidity.} Let
$V/{\mathbb P}^1$ be a Fano fiber space, satisfying the following
assumptions:

(i)  $V$ is a smooth projective variety with the Picard group
$\mathop{\rm Pic}V={\mathbb Z}K_V\oplus{\mathbb Z}F$, where $F$
is the class of a fiber of the projection $\pi\colon V\to{\mathbb
P}^1$,

(ii) every fiber $F_t=\pi^{-1}(t)$, $t\in{\mathbb P}^1$,  is a
(factorial) variety of dimension $\mathop{\rm dim}F_t\geq 4$,
with, at most, non-degenerate quadratic singularities, and
moreover,
$$
A^1F_t={\mathbb Z}H_t,\,\,A^2F_t={\mathbb Z}H^2_t,
$$
where $H_t=-K_{F_t}=(-K_V\cdot F_t)$ is the ample anticanonical
divisor,

(iii) at each point $x\in V$ the following local condition ($LR$)
holds: if $x\in F=F_t$ is a smooth point of the fiber, then for
every effective divisor $D\in|-nK_F|$ and every hyperplane
$B\subset E$ in the exceptional divisor $E\subset\widetilde F$ of
the blow up of the point $x$ on $F$ the following inequality is
satisfied:
\begin{equation}\label{c1a}
\mathop{\rm mult}\nolimits_xD+\mathop{\rm
mult}\nolimits_B\widetilde D\leq 2n,
\end{equation}
where $\widetilde D\subset\widetilde F$ is the strict transform
of the divisor $D$. If $x\in F$ is a quadratic singularity, then
for every effective divisor $D\in|-nK_F|$ and every hyperplane
section $B$ of the non-singular quadric $E\subset\widetilde F$,
the exceptional divisor of the blow up of the point $x$ on $F$,
the following inequality is satisfied:
\begin{equation}\label{c1}
\mathop{\rm mult}\nolimits_xD+2\mathop{\rm
mult}\nolimits_B\widetilde D\leq 4n,
\end{equation}
where $\widetilde D\subset\widetilde F$ is the strict transform
of the divisor $D$.

{\bf Proposition 3.1.} {\it For any movable linear system
$\Sigma\subset|-nK_V+lF|$ with $n\geq 1$ the log pair
$(V,\frac{1}{n}\Sigma)$ is canonical. In particular, if
$l\in{\mathbb Z}_+$, then the virtual and actual thresholds of
canonical adjunction coincide:}
$$
c_{\rm virt}(\Sigma)=c(\Sigma)=n.
$$

{\bf Proof.} Assume the converse: the pair
$(V,\frac{1}{n}\Sigma)$ is not canonical. In other words, the
linear system $\Sigma$ has a maximal singularity
$E\subset\widetilde V$, where $\varphi\colon\widetilde V\to V$ is
a sequence of blow ups. Let $Z=\varphi(E)\subset V$ be the centre
of the maximal singularity.

{\bf Lemma 3.1.} {\it The following estimate holds:} $\mathop{\rm
codim}_VZ\geq 3$.

{\bf Proof:} the claim follows immediately from the condition
(ii).

Let $x\in Z$ be a point of general position, $D\in\Sigma$ be a
general divisor. Furthermore, let
$$
\lambda\colon V^+\to V,
$$
be the blow up of the point $x$ on $V$. Denote the exceptional
divisor $\lambda^{-1}(x)\subset V^+$ by $E^+$.

{\bf Lemma 3.2.} {\it There exists a hyperplane $B^+\subset E^+$
satisfying the estimate
\begin{equation}\label{c2}
\mathop{\rm mult}\nolimits_xD+\mathop{\rm mult}\nolimits_{B^+}
D^+ > 2n,
\end{equation}
where $D^+\subset V^+$ is the strict transform of the divisor
$D$.}

{\bf Proof:} this follows easily from the connectedness principle
of Shokurov and Koll\' ar [4] and the properties of a
non-log-canonical isolated singularity on a surface. For a
detailed proof, see [3].

Let $F=F_t$ be the fiber of the projection $\pi$, containing the
point $x$. Assume first that $x\in F$ is a non-singular point of
the fiber. There are two possibilities:

1) $B^+\not\subset F^+$,

2) $B^+\subset F^+$,

\noindent
where $F^+\subset V^+$ is the strict transform of the
fiber $F^+$. Obviously, $\lambda_F\colon F^+\to F$ is the blow up
of the point $x\in F$, whereas $E_F=E^+\cap F^+$ is its
exceptional divisor. Consider these two possible cases separately.

1) Here $B=B^+\cap F^+$ is a hyperplane in the exceptional divisor
$E_F$. Set $D_F=D\,|\,_F$ to be the restriction of the divisor
$D$ on the fiber $F$. The effective divisor $D_F$ is well defined
since $F\not\subset\mathop{\rm Supp}D$. The inequality (\ref{c2})
can be rewritten in the following form:
$$
\mathop{\rm mult}\nolimits_{B^+}(\lambda^*D)> 2n
$$
(since, obviously, $\mathop{\rm mult}_{B^+}E^+=1$). Now we get
$$
\mathop{\rm mult}\nolimits_B(\lambda^*_FD_F)=\mathop{\rm
mult}\nolimits_B(\lambda^*D)\,|\,_{F^+}\geq\mathop{\rm
mult}\nolimits_{B^+}(\lambda^*D).
$$
Taking into account that $\mathop{\rm
mult}_B(\lambda^*_FD_F)=\mathop{\rm mult}_xD_F+\mathop{\rm
mult}_BD^+_F$, where $D^+_F\subset F^+$ is the strict transform
of the divisor $D_F$, we get finally
$$
\mathop{\rm mult}\nolimits_xD_F+\mathop{\rm
mult}\nolimits_BD^+_F>2n.
$$
Note that $D_F\in|-nK_F|$. We get a contradiction with the
condition ($LR$).

2) Here $B^+=E_F$ is the exceptional divisor of the blow up
$\lambda_F$. Setting $D_F=D\,|\,_F$, we obtain the estimate
$\mathop{\rm mult}\nolimits_xD_F=\mathop{\rm
mult}\nolimits_xD+\mathop{\rm mult}\nolimits_{B^+}D^+>2n$, which
contradicts the condition ($LR$) again.

Therefore the point $x\in F$ cannot be smooth. The essence of the
arguments above is that the inequality (\ref{c2}) can be
``restricted onto the fiber $F$''.

However, the point $x\in F$ cannot be a singularity of the fiber,
either. Indeed, in the notations above, in this case $E_F\subset
E^+$ is a non-singular quadric, so that $B^+\not\subset F^+$. The
argument 1) again yields the inequality
$$
\mathop{\rm mult}\nolimits_B(\lambda^*_FD_F)>2n,
$$
where $B=B^+\cap F^+$ is a hyperplane section of the quadric
$E_F$. Taking into account that
$$
\mathop{\rm mult}\nolimits_B(\lambda^*_FD_F)=\frac12\mathop{\rm
mult}\nolimits_xD_F+\mathop{\rm mult}\nolimits_BD^+_F,
$$
we get a contradiction with the condition ($LR$) once again.

Therefore, the log pair $(V,\frac{1}{n}\Sigma)$ is canonical.

If $l\in{\mathbb Z}_+$, then $c(\Sigma)=n$. The standard
arguments of [7] show that the strict inequality
$$
c_{\rm virt}(\Sigma)<n
$$
implies non-canonicity of the pair $(V,\frac{1}{n}\Sigma)$.
Therefore, we get the equality $c_{\rm virt}(\Sigma)=n$.

Q.E.D. for Proposition 3.1.

{\bf Corollary 3.1.} {\it Assume that the Fano fiber space
$V/{\mathbb P}^1$ satisfies the conditions {\rm (i)-(iii)} and the
$K$-condition, that is,
$$
-K\not\in \mathop{\rm Int} A^1_{\rm mov}V.
$$
Then the fiber space $V/{\mathbb P}^1$ is birationally
superrigid.} \vspace{0.3cm}

%%%%%%%%%%%%%%%%%%%%%%%%%%%%%%%%%%%%%%%%%%%%%%%%%%%%%%%%%%%%%%%%%%%
%%%%%%%%%%%%%%%%%%% subsection 3.2

{\bf 3.2. Varieties with a pencil of double spaces.} In this
section of the paper the symbol ${\mathbb P}$ stands for the
complex projective space ${\mathbb P}^M$ of dimension $M\geq 5$.
Let ${\cal W}={\mathbb P}(H^0({\mathbb P},{\cal O}_{\mathbb
P}(2M)))$ be the space of hypersurfaces of degree $2M$. In [3]
the following fact was proved:

{\bf Proposition 3.2.} {\it There exists a non-empty Zariski open
subset ${\cal W}_{\rm reg}\subset{\cal W}$ such that:

{\rm (i)} $\mathop{\rm codim}({\cal W}\setminus{\cal W}_{\rm
reg})\geq 2$,

{\rm (ii)} the Fano double space $\sigma\colon F\to{\mathbb P}$,
branched over a hypersurface $W\in{\cal W}_{\rm reg}$, satisfies
the condition ($LR$) at every smooth point $x\in F$ and has at
most non-denerate quadratic singularities.}

Assume that $x\in F$ is a singular point.

{\bf Lemma 3.3.} {\it The variety $F$ satisfies the condition
($LR$) at the point $x$.}

{\bf Proof.} This is almost obvious: if for a divisor
$D\in|-nK_F|$ and a hyperplane section $B$ of the quadric
$E=\lambda^{-1}(x)$, where $\lambda\colon\widetilde F\to F$ is
the blow up of the point $x$, the inequality
\begin{equation}\label{c4}
\mathop{\rm mult}\nolimits_xD+2\mathop{\rm
mult}\nolimits_B\widetilde D>4n,
\end{equation}
holds, then $\mathop{\rm mult}_xD>2n$, since $\mathop{\rm
mult}_xD\geq 2\mathop{\rm mult}_B\widetilde D$. However, the
opposite estimate holds: $\mathop{\rm mult}_xD\leq 2n$. This
contradiction completes the proof of the lemma.

Since $\mathop{\rm codim}({\cal W}\setminus {\cal W}_{\rm
reg})\geq 2$, for a general variety $V/{\mathbb P}^1$ with a
pencil of Fano double spaces of dimension $M\geq 5$ every fiber
satisfies the condition ($LR$) at every point.

Recall that the variety $V$ is the double cover of the projective
bundle $X/{\mathbb P}^1$, $X={\mathbb P}({\cal E})$, branched
over a smooth hypersurface $W^*\sim 2ML_X+2a_WR$, where $L_X$,
$R$ are the classes of the tautological sheaf and the fiber in
$\mathop{\rm Pic}X$, respectively. By construction, the locally
free sheaf ${\cal E}$ is generated by its sections, so that the
linear system $|L_V|$ is free, where $L_V=\sigma^*L_X$.

{\bf Lemma 3.4.} {\it The following equality holds: $(-K_V\cdot
L^{M-1}_V)=4-2a_W$.}

{\bf Proof:} direct computations.

{\bf Corollary 3.2.} {\it For $a_W\geq 2$ the Fano fiber space
$V/{\mathbb P}^1$ satisfies the $K$-condition:}
$-K_V\not\in\mathop{\rm Int}A^1_{\rm mov}V$.

In fact, for $a_W\geq 2$ a stronger condition holds:
$$
-K_V\not\in\mathop{\rm Int}A^1_+V,
$$
but we do not need that.

{\bf Corollary 3.3.} {\it For $a_W\geq 2$ a general Fano fiber
space $V/{\mathbb P}^1$ is birationally superrigid.}
\vspace{0.3cm}

%%%%%%%%%%%%%%%%%%%%%%%%%%%%%%%%%%%%%%%%%%%%%%%%%%%%%%%%%%%%%
%%%%%%%%%%%%%%%%%% subsection 3.3

{\bf 3.3. Varieties with a pencil of double quadrics.} Now set
${\mathbb P}={\mathbb P}^{M+1}$, $M\geq 6$. Set also ${\cal
W}={\mathbb P}(H^0({\mathbb P}, {\cal O}_{\mathbb P}(2M-2)))$. For
a hypersurface $W\in{\cal W}$ consider the double cover of a
quadric $Q\subset{\mathbb P}$:
$$
\sigma=\sigma_W\colon F\to Q,
$$
branched over the intersection $W_Q=W\cap Q$. About the quadric
$Q$, the hypersurface $W$ and their position with respect to each
other we assume the following.

1) The quadric $Q$ is either smooth or a cone over a smooth
quadric in ${\mathbb P}^M$ with the vertex at a point.

2) The variety $F$ has at most non-degenerate quadratic
singularities. If $Q$ is a cone, then the branch divisor $W_Q$ is
smooth.

3) Let $p\in W_Q$ be a smooth point of the branch divisor,
$(z_1,\dots,z_{M+1})$ a system of affine coordinates on ${\mathbb
P}^{M+1}$ with the origin at the point $p$. With respect to this
coordinate system the hypersurfaces $Q$ and $W$ are given by the
equations
\begin{equation}\label{c3}
q_1+q_2=0\,\, \mbox{and}\,\,w_1+w_2+\dots+w_{2M-2}=0,
\end{equation}
respectively. Without loss of generality assume that $q_1\equiv
z_1$, $w_1\equiv z_2$. Then for any linear form
$\lambda(z_3,\dots,z_{M+1})$ and any constant $c\in{\mathbb C}$
$$
c\,q_2(0,0,z_3,\dots,z_{M+1})\neq
w_2(0,0,z_3,\dots,z_{M+1})-\lambda^2(z_3,\dots,z_{M+1}).
$$
Similarly, let $p\in W_Q$ be a singularity of the branch divisor,
$(z_*)$ a system of affine coordinates, (\ref{c3}) the equations
of the hypersurfaces $Q$ and $W$, where without loss of
generality we assume that $q_1\equiv w_1\equiv z_1$. Then for any
linear form $\lambda(z_2,\dots,z_{M+1})$ and any constant
$c\in{\mathbb C}$
$$
c\,q_2(0,z_2,\dots,z_{M+1})\neq
w_2(0,z_2,\dots,z_{M+1})-\lambda^2(z_2,\dots,z_{M+1}).
$$

Let ${\cal Q}={\mathbb P}(H^0({\mathbb P},{\cal O}_{\mathbb
P}(2)))$ be the space of quadrics.

{\bf Proposition 3.3.} {\it There exists a Zariski open subset
${\cal F}_{\rm reg}\subset{\cal Q}\times{\cal W}$, such that

{\rm (i)} $\mathop{\rm codim}({\cal Q}\times{\cal
W}\setminus{\cal F}_{\rm reg})\geq 2$,

{\rm (ii)} for any pair $(Q,W)\in{\cal F}_{\rm reg}$ the
corresponding Fano double cover $F$ satisfies the conditions}
1)-3).

{\bf Proof:} an obvious dimension count.

{\bf Proposition 3.4.} {\it For any pair $(Q,W)\in{\cal F}_{\rm
reg}$ the corresponding Fano double cover $F$ satisfies the
condition ($LR$) at every point $x\in F$.}

{\bf Proof.} Set $p=\sigma(x)\in Q$. Consider first the case when
the morphism $\sigma$ is non-ramified at the point $x$, that is,
$p\not\in W$. Let
$$
\varphi\colon\widetilde F\to F \quad \mbox{and} \quad
\tilde\varphi\colon\widetilde Q\to Q
$$
be the blow ups of the points $x$ and $p$, respectively,
$E=\varphi^{-1}(x)$ and $\bar E=\bar\varphi^{-1}(p)$ the
exceptional divisors. The morphism $\sigma$ extends to a rational
map $\tilde\sigma\colon\widetilde F\dashrightarrow\widetilde Q$,
regular in a neighborhood of $E$ and identifying the exceptional
divisors $E$ and $\bar E$.

Assume that the point $p\in Q$ is smooth. Assume, moreover, that
the divisor $D\in|-nK_F|$ satisfies the inequality $\mathop{\rm
mult}\nolimits_xD+\mathop{\rm mult}\nolimits_B\widetilde D>2n$,
where $B\subset E$ is a hyperplane. Denote by the symbol $\bar
B=\tilde\sigma(B)$ the corresponding hyperplane in $\bar
E={\mathbb P}(T_pQ)$.

Now our arguments are similar to the proof of Theorem 2 in [3]
(with simplifications). Obviously, there is a pencil $\Lambda_B$
of hyperplane sections $R\ni x$ of the quadric $Q$, such that
$$
\widetilde R\cap\bar E=\bar B,
$$
where $\widetilde R\subset\widetilde Q$ is the strict transform
of the divisor $R$. A general quadric $R\in\Lambda_B$ is
non-singular at the point $p$, so that $\sigma^{-1}(R)$ is an
irreducible variety, non-singular at the point $x$. The class of
a hyperplane section of the variety $\sigma^{-1}(R)$ denote by
the symbol $H_R$. Set also $D_R=D\,|\,_{\sigma^{-1}(R)}$. This is
an effective divisor on the irreducible variety $\sigma^{-1}(R)$,
satisfying the inequality
$$
\mathop{\rm mult}\nolimits_x D_R=\mathop{\rm
mult}\nolimits_xD+\mathop{\rm mult}\nolimits_B\widetilde D>2n.
$$
However, the divisor $T_R=\sigma^{-1}(R\cap T_pR)$ is irreducible
and its multiplicity at the point $x$ is exactly 2. Therefore,
one can write down
$$
D_R=aT_R+D^{\sharp}_R,
$$
where $a\in{\mathbb Z}_+$ and the effective divisor
$D^{\sharp}_R\in|n^{\sharp}H_R|$ does not contain the divisor
$T_R$ as a component, and moreover $n^{\sharp}\geq 1$. We get the
inequality
$$
\mathop{\rm mult}\nolimits_xD^{\sharp}_R>2n^{\sharp}.
$$
Consider the effective cycle $\Delta=(D^{\sharp}_R\circ T_R)$.
Its degree with respect to the ample class $H_R$ is equal to
$4n^{\sharp}$, whereas the following inequality holds:
$$
\mathop{\rm mult}\nolimits_x\Delta>4n^{\sharp}.
$$
This is impossible. This contradiction excludes the case under
consideration.

Now assume that $p\in Q$ is the vertex of the cone. Here
$E\cong\bar E$ is a non-singular quadric, $B\cong\bar B$ is its
hyperplane section. Assume that some divisor $D\in|-nK_F|$
satisfies the inequality (\ref{c4}). There is a unique hyperplane
$H\subset{\mathbb P}$ such that the section $R=H\cap Q$ satisfies
the equality
$$
\widetilde R\cap\bar E=\bar B,
$$
$\widetilde R\subset\widetilde Q$ is the strict transform of the
section $R$. Since
$$
\mathop{\rm mult}\nolimits_xR=2,\,\,\mathop{\rm
mult}\nolimits_{\bar
B}R=1\,\,\mbox{and}\,\,\sigma^{-1}(R)\in|-K_F|,
$$
the divisor $\sigma^{-1}(R)\subset F$ satisfies the inequality
(\ref{c1}). Therefore, $\sigma^{-1}(R)\neq D$. Let us form the
effective cycle $\Delta=(D\circ\sigma^{-1}(R))$. Its
$(-K_F)$-degree is $4n$. Furthermore,
$$
\mathop{\rm mult}\nolimits_x\Delta=\mathop{\rm
mult}\nolimits_xD+2\mathop{\rm mult}\nolimits_B\widetilde D>4n.
$$
This is impossible. Thus the case $p\not\in W$ is completely
excluded.

Now assume that the morphism $\sigma$ is ramified at the point
$x$, that is, $p\in W$. The notations $\varphi\colon\widetilde
F\to F$, $\tilde\varphi\colon\widetilde Q\to Q$,
$E\subset\widetilde F$ and $\bar E\subset\widetilde Q$ mean the
same as above, but now the morphism $\sigma$  does {\it not}
identify $E$ and $\bar E$. Let $p\in Q$ be a smooth point,
$D\in|-nK_F|$ an effective divisor, satisfying the inequality
\begin{equation}\label{c5}
\mathop{\rm mult}\nolimits_xD+\mathop{\rm
mult}\nolimits_B\widetilde D>2n.
\end{equation}
Let $\bar\Lambda$ be the pencil of sections of the quadric $Q$ by
hyperplanes in ${\mathbb P}$, containing the linear space
$T_p(W\cap Q)=T_pW\cap T_p Q$. Take a general divisor $\bar
R\in\bar\Lambda$ and set $R=\sigma^{-1}(\bar R)$. Obviously,
$R\in|-K_F|$ is an irreducible divisor, $E_R=\widetilde R\cap E$
is an irreducible quadric in $E$. By the symbol $H_R$ we denote,
as above, the class of a hyperplane section of the variety $R$.
We get an effective divisor $D_R=(D\circ R)\in|nH_R|$ on $R$,
satisfying the equalities
$$
\mathop{\rm mult}\nolimits_xD_R=2\mathop{\rm mult}\nolimits_xD,
\quad \mathop{\rm mult}\nolimits_{B\cap E_R}\widetilde
D_R=\mathop{\rm mult}\nolimits_B\widetilde D.
$$
Since $\mathop{\rm mult}\nolimits_xD_R\leq 4n$, the multiplicity
$\mathop{\rm mult}\nolimits_B\widetilde D$ is strictly positive:
the support $\mathop{\rm Supp}\widetilde D$ contains the set $B$.
By linearity of the inequality (\ref{c5}) in the divisor $D$ we
can assume $D_R$ to be irreducible. As we have just explained,
$$
B\cap E_R\subset\mathop{\rm Supp}D_R.
$$

{\bf Lemma 3.5.} {\it The strict transform of the divisor
$T=\sigma^{-1}(\bar R\cap T_p\bar R)$ on $\widetilde F$ does not
contain the entire set $B\cap E_R$.}

{\bf Proof.} The claim follows from the conditions of general
position for the variety $F$. Indeed, set ${\mathbb T}={\mathbb
P}(T_p\bar R)$. It is easy to see that $\sigma\,|\,_R$ realizes
$E_R$ as the double cover of the space ${\mathbb T}$, branched
over the quadric
$$
W_T={\mathbb P}T_p(W\cap T_p(W\cap Q))
$$
(one should take into account that $T_p\bar R=T_p(W\cap Q))$. In
the coordinate form (with respect to the coordinate system in the
condition 3)) we get: the quadric $E_R$ is given by the equation
$$
y^2=w_2(0,0,z_3,\dots,z_{M+1})
$$
in the projective space with the homogeneous coordinates
$$
(y: z_3:\dots: z_{M+1}),
$$
$z_3:\dots: z_{M+1}$ gives a system of homogeneous coordinates on
${\mathbb T}$, the covering morphism $\sigma_E\colon
E_R\to{\mathbb T}$ is the restriction onto $E_R$ of the linear
projection from the point $(1: 0:\dots: 0)$. The set $\widetilde
T\cap E_R$ (where $\widetilde T$, as usual, is the strict
transform of the divisor $T$ on $\widetilde R$) is given by the
equation $q_2(0,0,z_3,\dots,z_{M+1})=0$.

Let $\alpha y+\lambda (z_3,\dots,z_{M+1})=0$ be the equation of
the hyperplane section $B\cap E_R$, $\alpha=0$ or 1,
$\lambda(\cdot)$ a linear form. Assume that
$$
B\cap E_R\subset\widetilde T.
$$
Then, obviously, $\sigma_E(B\cap E_R)$ is contained in the quadric
$q_2\,|\,_{\mathbb T}=0$. If $\alpha=0$, then it means that the
quadratic polynomial $q_2(0,0,z_*)$ is divisible by the linear
form $\lambda(z_*)$, which is impossible. If $\alpha=1$, then it
means that the quadratic polynomial $q_2(0,0,z_*)$ up to a scalar
factor is $w_2(0,0,z_*)-\lambda^2(z_*)$, which is again
impossible by the conditions of general position. Q.E.D. for the
lemma.

Since, as we have just shown, $D_R\neq T$ (and both divisors are
irreducible), the effective cycle $(D_R\circ T)=\Delta$ of
codimension 2 on $R$ is well defined. It is easy to see that its
degree with respect to $H_R$ is $4n$, whereas its multiplicity at
the point $x$ satisfies the inequality
$$
\mathop{\rm mult}\nolimits_x\Delta\geq 2\mathop{\rm
mult}\nolimits_xD_R=4\mathop{\rm mult}\nolimits_xD>4n,
$$
which is impossible. Thus the case of a smooth point $x\in F$,
$p\in W$ is excluded.

It remains to consider the only case of the singular point $x\in
F$, where $p=\sigma(x)\in W$. By the conditions of general
position we get: $p\in Q$ is a smooth point on the quadric.
Assume that a divisor $D\in|-nK_F|$ satisfies the inequality
(\ref{c4}). By linearity of this inequality in the divisor $D$ we
may assume that $D$ is a prime divisor. In the anticanonical
system $|-K_F|$ consider the divisor $T=\sigma^{-1}(Q\cap
T_{p}Q)$. Obviously, $\mathop{\rm mult}_xT=4$ and by the
conditions of general position $\widetilde T\not\supset B$.
Therefore, $D\neq T$ and the effective cycle $\Delta=(D\circ T)$
is well defined. Since
$$
\mathop{\rm mult}\nolimits_x\Delta\geq 2\mathop{\rm
mult}\nolimits_xD>4n
$$
($x\in F$ is a double point), we get a contradiction again. This
completes the proof of Proposition 3.4. Propositions 3.3, 3.4
imply that for a general Fano fiber space $V/{\mathbb P}^1$ into
double quadrics of index 1 and dimension $M\geq 6$ every fiber at
every point satisfies the condition ($LR$).

Recall [1,2], that the variety $V$ is realized as the double cover
$\sigma\colon V\to Q^*$ of the hypersurface $Q^*\subset X={\mathbb
P}({\cal E})$, branched over the divisor $W^*\cap Q^*$. The fiber
space $V$ is characterized by the triple of non-negative integral
parameters $(a_X,a_Q,a_W)$. It follows from the formula
$$
(-K_V\cdot L^M_V)=2(8-2a_X-3a_Q-4a_W)
$$
that if $2a_X+3a_Q+4a_W\geq 8$, then the Fano fiber space
$V/{\mathbb P}^1$ satisfies the $K$-condition:
$-K\not\in\mathop{\rm Int}A^1_{\rm mov}V$ (and even
$-K\not\in\mathop{\rm Int}A^1_+V$).

{\bf Corollary 3.4.} {\it For $2a_X+3a_Q+4a_W\geq 8$ a general
Fano fiber space $V/{\mathbb P}^1$ is birationally superrigid.}
\vspace{0.3cm}

%%%%%%%%%%%%%%%%%%%%%%%%%%%%%%%%%%%%%%%%%%%%%%%%%%%%%%%%%%%%%%%%%%
%%%%%%%%%%%%%%%%%%%%%%% subsection 3.4

{\bf 3.4. Fano double hypersurfaces.} We give here one
intermediate result for this class of varieties. Let
$\sigma\colon F\to Q\subset{\mathbb P}={\mathbb P}^{M+1}$ be a
Fano double cover, where $Q$ is a hypersurface of degree $m\geq
3$, $M\geq 6$. The morphism $\sigma$ is branched over the divisor
$W=W^*\cap Q$, where $W^*\subset{\mathbb P}$ is a hypersurface of
degree $2l$, $m+l=M+1$. Assume that $F$ has at most
non-degenerate quadratic singularities. Assume also that at any
smooth point $x\in F$ outside the ramification divisor the variety
$F$ is regular in the sense of Definition 1 of the paper [9].
Besides, assume that the following conditions of general position
are satisfied. Let $(z_1,\dots,z_{M+1})$ be a system of affine
coordinates with the origin at the point $p=\sigma(x)$,
$$
f=q_1+q_2+\dots+q_m
$$
an equation of the hypersurface $Q$ (so that by Definition 1 of
the paper [9] the sequence $q_1,\dots,q_m$ is regular). We
require that the linear span of any irreducible component of the
closed set
$$
q_1=q_2=q_3=0
$$
in ${\mathbb C}^{M+1}$ coincides with the hyperplane $q_1=0$
(that is, the tangent hyperplane $T_{p}Q$). Furthermore, we
require that the closed algebraic set
$$
\sigma^{-1}(\overline{\{q_1=q_2=0\}\cap Q})\subset F
$$
should be irreducible and any section of this set by an
anticanonical divisor $P\in|-K_F|$, $P\ni x$, should be also
irreducible and reduced.

These conditions are completely similar to the regularity
conditions (R1.1)-(R1.3) of the paper [3].

{\bf Proposition 3.5.} {\it In these assumptions the variety $F$
satisfies the condition ($LR$) at any smooth point $x\in F$
outside the ramification divisor.}

To obtain the {\bf proof}, we repeat the arguments in [3, Sec.
2.1] word for word, replacing (in the concluding part of the
proof) the hypertangent technique of the paper [5] (which is used
in [3]) by the hypertangent technique of the paper [9]. No other
changes are necessary. Q.E.D. for the proposition.

{\bf Corollary 3.5.} {\it Assume that the Fano fiber space
$V/{\mathbb P}^1$ into double hypersurfaces of index 1 satisfies
the conditions of general position, formulated above, at every
smooth point of every fiber, lying outside the ramification
divisor. Then the centre $B=\mathop{\rm centre}(E)$ of any
maximal singularity $E$ of a movable linear system
$\Sigma\subset|-nK_V+lF|$ with $l\in{\mathbb Z}_+$ either is
contained in the ramification divisor, or coincides with a
singular point of a fiber.}\vspace{0.3cm}

%%%%%%%%%%%%%%%%%%%%%%%%%%%%%%%%%%%%%%%%%%%%%%%%%%%%%%%%%%%%%%
%%%%%%%%%%%%%%%%%%%% subsection 3.5

{\bf 3.5. Varieties with a pencil of Fano hypersurfaces.} To
conclude this section, let us consider the class of Fano fiber
spaces, that were studied in [7]. Let us show how, at the expense
of making the conditions of general position stronger, we can use
the linear method to simplify the proof of birational rigidity.

Let ${\cal F}={\mathbb P}(H^0({\mathbb P}, {\cal O}_{\mathbb
P}(M)))$, ${\mathbb P}={\mathbb P}^M$, be the set of Fano
hypersurfaces of index 1. In [3, Sec. 2] the following fact was
obtained.

{\bf Proposition 3.6.} {\it There exists an open subset (in the
sense of Zariski topology) ${\cal F}_{\rm reg}\subset{\cal F}$
such that any hypersurface $F\in{\cal F}_{\rm reg}$ is smooth and
satisfies the condition ($LR$) at every point.}

Singular hypersurfaces form a closed subset in ${\cal F}$ of
codimension 1, so that the complement ${\cal F}\setminus{\cal
F}_{\rm reg}$ is of codimension 1, too. However, from the
computations of [3, Sec. 2.3] one can see that violation of the
regularity conditions $(R1.1)-(R1.3)$ (which define the set ${\cal
F}_{\rm reg}$) at at least one smooth point $x\in F$ defines a
closed subset of ${\cal F}$ of codimension at least 2. Therefore,
the following fact takes place.

{\bf Proposition 3.7.} {\it There exists an open (in the sense of
Zariski topology) subset ${\cal F}^+_{\rm reg}\subset{\cal F}$
such that

{\rm (i)} $\mathop{\rm codim}({\cal F}\setminus{\cal F}^+_{\rm
reg})\geq 2$,

{\rm (ii)} any hypersurface ${\cal F}\in{\cal F}^+_{\rm reg}$ has
at most non-degenerate quadratic singularities, and moreover, if
$x\in F$ is such a point, then the sequence of homogeneous
polynomials
$$
q_2,\dots,q_M
$$
is regular in ${\cal O}_{x,{\mathbb P}}$, where $f=q_2+\dots+q_M$
is the polynomial, defining the hypersurface $F$ in a system of
affine coordinates with the origin at the point $x$,

{\rm (iii)} any hypersurface ${\cal F}\in{\cal F}^+_{\rm reg}$
satisfies the condition ($LR$) at every smooth point $x\in F$.}

It turns out that non-degenerate quadratic singularities,
satisfying the regularity condition (ii), generate no serious
problems.

{\bf Proposition 3.8.} {\it Any hypersurface ${\cal F}\in{\cal
F}^+_{\rm reg}$ satisfies the condition ($LR$) at every point
$x\in F$.}

{\bf Proof.} We need to check only the case of a quadratic
singularity $x\in F$. Let
$$
\varphi\colon\widetilde{\mathbb P}\to {\mathbb
P}\,\,\mbox{and}\,\,\varphi_F\colon\widetilde F\to F
$$
be the blow up of the point $x\in{\mathbb P}$ and its restriction
onto $F$, respectively,
$$
E=\varphi^{-1}(x)\subset\widetilde{\mathbb
P}\,\,\mbox{and}\,\,E_F=\varphi^{-1}_F(x)=\widetilde F\cap E
$$
the exceptional divisors,  $E_F\subset E$ a smooth quadric.
Assume that the divisor $D\in|-nK_F|$ satisfies the inequality
(\ref{c4}). By linearity of this inequality in $D$ we may assume
that the divisor $D$ is irreducible. We get
$$
B=E_F\cap\Pi,
$$
where $\Pi\subset E$ is a hyperplane. Let $P\subset{\mathbb P}$,
$P\ni x$ be the only hyperplane, for which the equality
$\widetilde P\cap E=\Pi$ holds, $\widetilde
P\subset\widetilde{\mathbb P}$ is the strict transform. Set
$R=P\cap F$. This is a divisor on $F$.

Obviously, $R$ is irreducible and satisfies the inequality
(\ref{c1}) with $n=1$ and $R$ instead of $D$. Therefore, $R\neq
D$, so that the effective cycle $Y=(R\circ D)$ of codimension 2
is well defined. This cycle (we may assume it to be irreducible)
satisfies the inequality
$$
\frac{\mathop{\rm mult}_x}{\mathop{\rm deg}}Y>\frac{4}{M}.
$$
Now the proof is completed in the standard way: let
$\Lambda_2,\dots,\Lambda_{M-1}$ be the hypertangent linear
systems. By the regularity condition we get
$$
\mathop{\rm codim}\nolimits_F\mathop{\rm Bs}\Lambda_k=k-1.
$$
Let
$(D_4,\dots,D_{M-1})\in\Lambda_4\times\dots\times\Lambda_{M-1}$
be a generic set of hypertangent divisors. Obviously,
$Y\not\subset\mathop{\rm Supp}\Lambda_4$. Let us construct by
induction a sequence of irreducible subvarieties $Y_i$,
$i=2,\dots,M-2$, where $\mathop{\rm codim}_FY_i=i$, $Y_2=Y$,
satisfying the inequalities
$$
\frac{\mathop{\rm mult}_x}{\mathop{\rm
deg}}Y_{i+1}\geq\frac{\mathop{\rm mult}_x}{\mathop{\rm
deg}}Y_i\cdot\frac{i+3}{i+2},
$$
$Y_{i+1}$ is an irreducible component of the effective cycle
$(Y_i\circ D_{i+2})$. This is possible, because
$$
\mathop{\rm codim}\nolimits_FY_i=i<\mathop{\rm
codim}\nolimits_F\mathop{\rm Bs}\Lambda_{i+2}=i+1,
$$
so that $Y_i\not\subset\mathop{\rm Supp}D_{i+2}$ for a general
divisor $D_{i+2}\in\Lambda_{i+2}$. The curve $C=Y_{M-2}$
satisfies the inequality
$$
\frac{\mathop{\rm mult}_x}{\mathop{\rm
deg}}C>\frac{4}{M}\cdot\frac{5}{4}\cdot\dots\cdot\frac{M}{M-1}=1,
$$
which is impossible. The contradiction obtained above completes
the proof of Proposition 3.8.

{\bf Corollary 3.6.} {\it Assume that for every fiber
$F_t=\pi^{-1}(t)$ of the Fano fiber space $\pi\colon V\to{\mathbb
P}^1$ we have $F_t\in{\cal F}^+_{\rm reg}$. Then:

{\rm (i)} for any movable linear system $\Sigma\subset|-nK_V+lF|$
with $l\in{\mathbb Z}_+$ the pair $(V,\frac{1}{n}\Sigma)$ is
canonical,

{\rm (ii)} if the variety $V$ satisfies the $K$-condition, then
it is birationally superrigid.}

{\bf Remark 3.1.} The linear method makes it possible to simplify
the proofs of birational rigidity in [7,10]. However, the
technique used for proving the conditions of the type ($L$),
works for varieties of sufficiently high (at least 3) dimension.
For this reason, the linear method in its present form gives not
much for fibrations into curves [11,12] and surfaces [13-15]. For
instance, in the case of one-dimensional fibers, restricting the
linear system onto a conic, intersecting the centre of a maximal
singularity, we get no new information. In the case of
two-dimensional fibers (del Pezzo surfaces) the situation is
better, however, it is easy to show that under any conditions of
general position there are a fiber $F$, a point $x\in F$ and a
divisor $D\in|-nK_F|$ such that the pair $(F,\frac{1}{n}D)$ is
not log canonical at the point $x$. Thus it is impossible to
exclude a maximal singularity with the centre at this point by
the linear method. That the linear method works in smaller
dimensions with difficulties, can be seen already in the absolute
case (that is, for Fano varieties): the condition ($L$) remains
unproved for three-dimensional quartics and four-dimensional
quintics, whereas the quadratic method works quite successfully
for these varieties [16,17].

%%%%%%%%%%%%%%%%%%%%%%%%%%%%%%%%%%%%%%%%%%%%%%%%%%%%%%%%%%%%%%%%%%
%%%%%%%%%%%%%%%%%%%%%%%%%%%%%%%%%%%%%%%%%%%%%%%%%%%%%%%%%%%%%%%%%%
%%%%%%%%%%%%%%%%%%%%%%%%%%%%%%%%%%%%%%%%%%%%%%%%%%%%%%%%%%%%%%%%%%

%%%%%%%%%%%%%%%%%%%%%%%%%%%%%%%%%%%%%%%%%%%%%%%%%%%%%%%%%%%%%%%%%%%
%%%%%%%%%%%%%%%%%%%%%%%%%%%%%%%%%%%%%%%%%%%%%%%%%%%%%%%%%%%%%%%%%%%
%%%%%%%%%%%%%%%%%%%%%%%%%%%%%%%%%%%%%%%%%%%%%%%%%%%%%%%%%%%%%%%%%%%
%%%%%%%%%%%%%%%%%%%%%%%%%%%%%%%%%%%%%%%%%%%%%%%%%%%%%%%%%%%%%%%%%%%
%%%%%%%%%%%%%%%%%%    references                       %%%%%%%%%%%%
\section*{References}
{\small

\noindent 1. Pukhlikov A.V., Birationally rigid varieties with a
pencil of Fano double covers. I. Sbornik: Mathematics {\bf 195}
(2004), no. 7, 1039-1071. \vspace{0.3cm}

\noindent 2. Pukhlikov A.V., Birationally rigid varieties with a
pencil of Fano double covers. II. Sbornik: Mathematics {\bf 195}
(2004), no. 11, 1665-1702. \vspace{0.3cm}

\noindent 3. Pukhlikov A.V., Birational geometry of Fano direct
products, arXiv: math.AG/0405011.\vspace{0.3cm}

\noindent 4. Koll{\'a}r J., et al., Flips and Abundance for
Algebraic Threefolds, Asterisque 211, 1993. \vspace{0.3cm}

\noindent 5. Pukhlikov A.V., Birational automorphisms of Fano
hypersurfaces, Invent. Math. {\bf 134} (1998), no. 2, 401-426.
\vspace{0.3cm}

\noindent 6. Pukhlikov A.V., Fiber-wise birational
correspondences, Mathematical Notes {\bf 68} (2000), no. 1,
103-112. \vspace{0.3cm}

\noindent 7. Pukhlikov A.V., Birationally rigid Fano fibrations,
Izvestiya: Mathematics {\bf 64} (2000), 131-150. \vspace{0.3cm}

\noindent 8. Fulton W., Intersection Theory, Springer-Verlag,
1984. \vspace{0.3cm}

\noindent 9. Pukhlikov A.V., Birationally rigid Fano double
hypersurfaces, Sbornik: Mathematics {\bf 191} (2000), No. 6,
101-126. \vspace{0.3cm}

\noindent 10. Sobolev I. V., On a series of birationally rigid
varieties with a pencil of Fano hypersurfaces, Sbornik:
Mathematics {\bf 192} (2001), no. 9-10, 1543-1551. \vspace{0.3cm}

\noindent 11. Sarkisov V.G., Birational automorphisms of conic
bundles, Math. USSR Izv. {\bf 17} (1981), 177-202. \vspace{0.3cm}

\noindent 12. Sarkisov V.G., On conic bundle structures,  Math.
USSR Izv. {\bf 20} (1982), no. 2, 354-390. \vspace{0.3cm}

\noindent 13. Pukhlikov A.V., Birational automorphisms of
three-dimensional algebraic varieties with a pencil of del Pezzo
surfaces, Izvestiya: Mathematics {\bf 62}:1 (1998), 115-155.
\vspace{0.3cm}

\noindent 14. Grinenko M.M., Birational properties of pencils of
del Pezzo surfaces of degrees 1 and 2. Sbornik: Mathematics. {\bf
191} (2000), no. 5, 17-38. \vspace{0.3cm}

\noindent 15. Grinenko M.M., Birational properties of pencils of
del Pezzo surfaces of degrees 1 and 2. II. Sbornik: Mathematics.
{\bf 194} (2003). \vspace{0.3cm}

\noindent 16. Iskovskikh V.A. and Manin Yu.I., Three-dimensional
quartics and counterexamples to the L\" uroth problem, Math. USSR
Sb. {\bf 86} (1971), no. 1, 140-166. \vspace{0.3cm}

\noindent 17. Pukhlikov A.V., Birational isomorphisms of
four-dimensional quintics, Invent. Math. {\bf 87} (1987), 303-329.
\vspace{0.3cm} }

\begin{flushleft}
{\it e-mail}: pukh@liv.ac.uk, pukh@mi.ras.ru
\end{flushleft}

\end{document}